\documentclass[a4paper]{amsart}
\usepackage{amsmath,amssymb}
\usepackage{graphicx}
\usepackage[all]{xy}

\usepackage{multirow}
\usepackage{color}
\newtheorem{Theorem}{Theorem}[section]
\newtheorem{Lemma}[Theorem]{Lemma}
\newtheorem{Proposition}[Theorem]{Proposition}

\newtheorem{Example}[Theorem]{Example}

\newtheorem{Conjecture}[Theorem]{Conjecture}

\newtheorem{Remark}[Theorem]{Remark}

\def\labelenumi{\rm ({\makebox[0.8em][c]{\theenumi}})}
\def\theenumi{\arabic{enumi}}
\def\labelenumii{\rm {\makebox[0.8em][c]{(\theenumii)}}}
\def\theenumii{\roman{enumii}}

\renewcommand{\thefigure}
{\arabic{section}.\arabic{figure}}
\makeatletter
\@addtoreset{figure}{section}
\def\@thmcountersep{-}
\makeatother

\numberwithin{equation}{section}

\usepackage{here}
\begin{document}

\title[Pairs of knot invariants]{Pairs of knot invariants}

\author{Kouki Taniyama}
\address{Department of Mathematics, School of Education, Waseda University, Nishi-Waseda 1-6-1, Shinjuku-ku, Tokyo, 169-8050, Japan}
\email{taniyama@waseda.jp}

\thanks{The author was partially supported by Grant-in-Aid for Scientific Research (c) (No. 21K03260) and Grant-in-Aid for Scientific Research (A) (No. 21H04428), Japan Society for the Promotion of Science.}


\date{}

\dedicatory{}


\begin{abstract}

Let $\alpha$ be a map from the set of all knot types ${\mathcal K}$ to a set $X$. Let $\beta$ be a map from ${\mathcal K}$ to a set $Y$. 
We define the relation between $\alpha$ and $\beta$ to be the image of a map $(\alpha,\beta)$ from ${\mathcal K}$ to $X\times Y$ sending an element $K$ of ${\mathcal K}$ to $(\alpha(K),\beta(K))$. 
We determine the relations between $\alpha$ and $\beta$ for certain $\alpha$ and $\beta$ such as crossing number, unknotting number, bridge number, braid index, genus and canonical genus. This is a study of geography problem in knot theory. 

\end{abstract}

\maketitle

\section{Introduction}\label{introduction} 

Throughout this paper we denote the set of all oriented knot types in the $3$-sphere ${\mathbb S}^{3}$ by ${\mathcal K}$. 
We do not distinguish a knot and its knot type so long as no confusion occurs. 
We regard a map from ${\mathcal K}$ to a set $X$ as a knot invariant taking its values in $X$.  
Let 
\[
c:{\mathcal K}\to{\mathbb Z}_{\geq0}
\]
be the crossing number and 
\[
u:{\mathcal K}\to{\mathbb Z}_{\geq0}
\]
the unknotting number. 
Namely, the crossing number of an oriented knot $K$ is denoted by $c(K)$ and the unknotting number of $K$ is denoted by $u(K)$. 
It is well-known that every nontrivial knot $K$ satisfies the inequality 
\[
\displaystyle{u(K)\leq\frac{1}{2}(c(K)-1)}.
\]
The equality holds if and only if $K$ is a $(2,n)$-torus knot for some odd number $n$ \cite{Taniyama}. Therefore this inequality is considered to be the only relation between crossing number and unknotting number of a knot. 
However, the trivial knot does not satisfy this inequality, and there may be some other entirely different inequalities between them. For example, there may be an inequality that estimates unknotting number below by crossing number. 
In fact, the following is an almost trivial but nontrivial estimate below for any knot $K$. 
\[
u(K)\geq{\rm min}\{c(K),1\}.
\]
In order to make it explicit, we define the relation between two knot invariants as follows. 
We define it in a more general setting. We define the relation between two topological invariants of some topological objects. 

Let ${\mathcal T}$ be a set. Let $X$ and $Y$ be sets and $\alpha:{\mathcal T}\to X$, $\beta:{\mathcal T}\to Y$ maps. 
Let 
\[
(\alpha,\beta):{\mathcal T}\to X\times Y
\]
be a map defined by $(\alpha,\beta)(T)=(\alpha(T),\beta(T))$ for $T\in{\mathcal T}$. We call $(\alpha,\beta)$ the {\it pair} of $\alpha$ and $\beta$. The subset $(\alpha,\beta)({\mathcal T})$ of $X\times Y$ is said to be the {\it relation} between $\alpha$ and $\beta$.

\begin{figure}[htbp]
      \begin{center}
\scalebox{0.6}{\includegraphics*{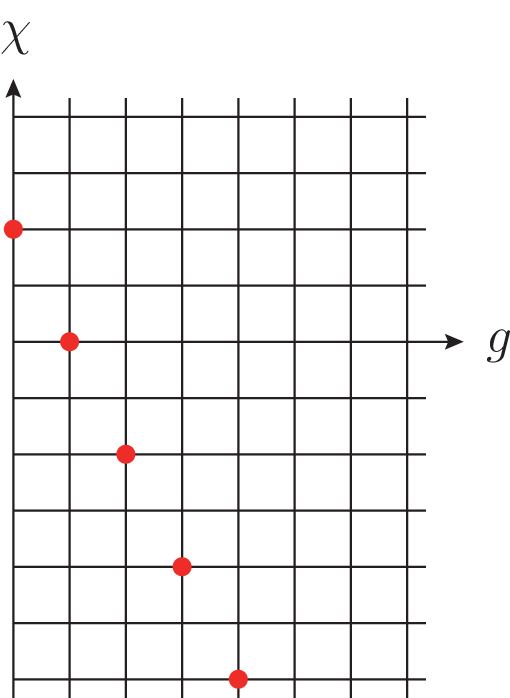}}
      \end{center}
   \caption{The relation between genus and Euler characteristic of closed orientable surfaces}
  \label{genusEuler}
\end{figure}

\begin{Example}\label{example-genus-Euler-characteristic}
{\rm
Let ${\mathcal T}$ be the set of all homeomorphism classes of closed connected orientable surfaces. 
Let $g:{\mathcal T}\to{\mathbb Z}_{\geq0}$ be the genus and $\chi:{\mathcal T}\to{\mathbb Z}$ the Euler characteristic. 
We note that $\chi(F)=2-2g(F)$ for each $F\in{\mathcal T}$. Therefore the relation between genus and Euler characteristic is the set 
\[
\{(x,y)\in{\mathbb Z}_{\geq0}\times{\mathbb Z}\mid y=2-2x\}.
\]
See Figure \ref{genusEuler}. 
}
\end{Example}

The study of $(\alpha,\beta)({\mathcal T})$ is called geography problem in some areas of mathematics. 
It begins with \cite{Persson} where the pair $(c_{1}^{2},c_{2})$ of Chern numbers of a surface of general type is considered. 
The study of $(\alpha,\beta)^{-1}(x,y)$ for an element $(x,y)$ of $X\times Y$ is called botany problem. 
In this paper we study geography problem for certain pairs of knot invariants.

From now on we restrict our attention to knot invariants. 
We note that the relation between two knot invariants defined above is a natural concept and implicitly used for a long time. 
Prior to this concept for a pair of knot invariants, a concept for a single knot invariant exists. 
Let $X$ be a set and $\alpha:{\mathcal K}\to X$ a knot invariant. 
Then the determination of the image $\alpha({\mathcal K})\subset X$ is a fundamental problem. 
In many cases it is easily determined. 
For example, let $\Delta:{\mathcal K}\to{\mathbb Z}[t^{\pm1}]$ be the Alexander polynomial. 
It is known that 
\[
\Delta({\mathcal K})=\{p(t)\in{\mathbb Z}[t^{\pm1}]\mid p(t^{-1})=p(t),p(1)=1\}. 
\]
This is known to be a characterization of Alexander polynomial. 
However, a similar characterization of, say Jones polynomial is not known yet. 

\begin{Example}\label{example-unknotting-number-Alexander-polynomial}
{\rm
It is independently shown in \cite{Kondo} and \cite{Sakai} that 
\[
\Delta({\mathcal K})=\Delta(\{K\in{\mathcal K}\mid u(K)=1\}).
\]
We note that this is equivalent to 
\[
(u,\Delta)({\mathcal K})\cap(\{1\}\times{\mathbb Z}[t^{\pm1}])=\{1\}\times\Delta({\mathcal K}). 
\]
Then we see that the Alexander polynomial $\Delta(K)$ of an oriented knot $K$ never tell us $u(K)\geq2$. 
Thus the set $(u,\Delta)({\mathcal K})$ has information that the lower estimate of unknotting number by Alexander polynomial alone is at most $1$. 
}
\end{Example}

There are a lot of knot invariants. In this paper we mainly consider crossing number, unknotting number, braid index, bridge number, genus and canonical genus. 
Most of the results in this paper are expected ones. Existence of certain elements of $(\alpha,\beta)({\mathcal K})$ is already known for $\alpha$ and $\beta$ one of above six knot invariants. However, as far as the author knows, $(\alpha,\beta)({\mathcal K})$ itself was not explicitly decided yet except the pair of canonical genus and genus that is described in \cite{Stoimenow} without proof. 

We use symbols in Rolfsen's knot table \cite{Rolfsen} and Hoste-Thistlethwaite's table of 11 Crossing Knots \cite{H-T}\cite{H-T-W}. 
In particular, $0_{1}$ is a trivial knot, $3_{1}$ is a left-handed trefoil knot and $4_{1}$ is a figure-eight knot. 
We denote the mirror image of a knot $K$ by $K^*$. 
A torus knot of type $(p,q)$ is denoted by $T(p,q)$ and a pretzel knot of type $(p_{1},\cdots,p_{k})$ is denoted by $P(p_{1},\cdots,p_{k})$. 
A $2$-bridge knot is denoted by Conway notation $C(p_{1},\cdots,p_{k})$ so that if all of $p_{1},\cdots,p_{k}$ are positive, then the diagram corresponding to this notation is alternating, and therefore it is a minimal crossing diagram. Namely we have $c(C(p_{1},\cdots,p_{k}))=p_{1}+\cdots+p_{k}$. 

We denote the connected sum of two oriented knots $J$ and $K$ by $J\#K$. 
We denote the connected sum of $p$ copies of $K$ by $p\cdot K$. 
We denote the number of elements of a finite set $X$ by $|X|$. 

Let ${\mathcal P}$ be the set of all oriented prime knot types and let ${\mathcal P}_{0}={\mathcal P}\cup\{0_{1}\}$. 
Let ${\mathcal A}$ be the set of all oriented alternating knot types. We note that $0_{1}$ is an element of ${\mathcal A}$. 
Let ${\mathcal R}$ be the set of all oriented $2$-bridge knot types and let ${\mathcal R}_{0}={\mathcal R}\cup\{0_{1}\}$. 
We note that both ${\mathcal P}_{0}$ and ${\mathcal A}$ are proper subsets of ${\mathcal K}$ and ${\mathcal R}_{0}$ is a proper subset of ${\mathcal P}_{0}\cap{\mathcal A}$. 

We denote the bridge number and braid index of an oriented knot $K$ by ${\rm bridge}(K)$ and ${\rm braid}(K)$ respectively. Namely we have knot invariants 
${\rm bridge}:{\mathcal K}\to{\mathbb Z}_{>0}$ and ${\rm braid}:{\mathcal K}\to{\mathbb Z}_{>0}$. We define knot invariants 
\[
({\rm bridge}-1):{\mathcal K}\to{\mathbb Z}_{\geq0}
\]
and 
\[
({\rm braid}-1):{\mathcal K}\to{\mathbb Z}_{\geq0}
\]
by $({\rm bridge}-1)(K)={\rm bridge}(K)-1$ and $({\rm braid}-1)(K)={\rm braid}(K)-1$ respectively as a slight modification of them. We call them {\it bridge number minus one} and {\it braid index minus one} respectively. Since 
\[
{\rm bridge}(J\#K)={\rm bridge}(J)+{\rm bridge}(K)-1
\]
\cite{Schubert} \cite{Schultens} and 
\[
{\rm braid}(J\#K)={\rm braid}(J)+{\rm braid}(K)-1
\]
\cite{B-M}, we have 
\[
({\rm bridge}-1)(J\#K)=({\rm bridge}-1)(J)+({\rm bridge}-1)(K)
\]
and 
\[
({\rm braid}-1)(J\#K)=({\rm braid}-1)(J)+({\rm braid}-1)(K).
\] 
Namely these modified invariants are additive under connected sum of knots.

Let 
\[
g:{\mathcal K}\to{\mathbb Z}_{\geq0}
\]
be the genus and 
\[
g_{c}:{\mathcal K}\to{\mathbb Z}_{\geq0}
\]
the canonical genus. 
Namely, the genus of an oriented knot $K$ is denoted by $g(K)$ and the canonical genus of $K$ is denoted by $g_{c}(K)$. 
By definition we have $g(K)\leq g_{c}(K)$ for any oriented knot $K$.

First we describe the relation between crossing number and one of above invariants. 
The relations between crossing number and unknotting number, crossing number and genus, and crossing number and canonical genus are all equal to a set described in the following theorem.

\begin{Theorem}\label{theorem-crossing-number-unknotting-number-genus-canonical-genus}
\[
(c,u)({\mathcal K})=(c,g)({\mathcal K})=(c,g_{c})({\mathcal K})=\{(0,0)\}\cup\{(x,y)\in({\mathbb Z}_{>0})^{2}\mid y\leq\frac{1}{2}(x-1)\}.
\]
Moreover
\[
(c,u)({\mathcal K})=\{(6,2)\}\cup(c,u)({\mathcal P}_{0})=\{(c,u)({3_{1}}^*\#{3_{1}}^*)\}\cup(c,u)({\mathcal P}_{0})
\]
and
\[
(c,g)({\mathcal K})=(c,g)({\mathcal R}_{0})=(c,g_{c})({\mathcal R}_{0}).
\]

\end{Theorem}

\vskip 5mm

\begin{figure}[htbp]
      \begin{center}
\scalebox{0.6}{\includegraphics*{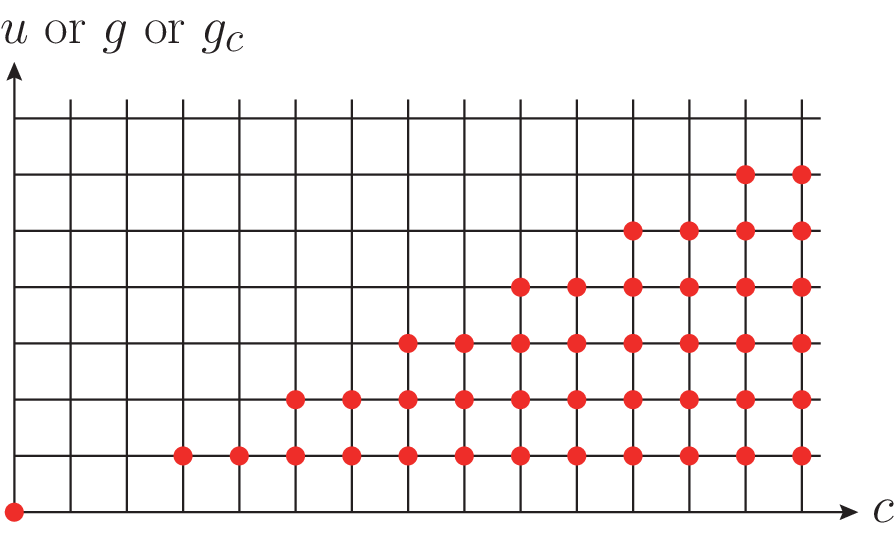}}
      \end{center}
   \caption{$(c,u)({\mathcal K})=(c,g)({\mathcal K})=(c,g_{c})({\mathcal K})$}
  \label{cugcg}
\end{figure}

See Figure \ref{cugcg}.

\begin{figure}[htbp]
      \begin{center}
\scalebox{0.4}{\includegraphics*{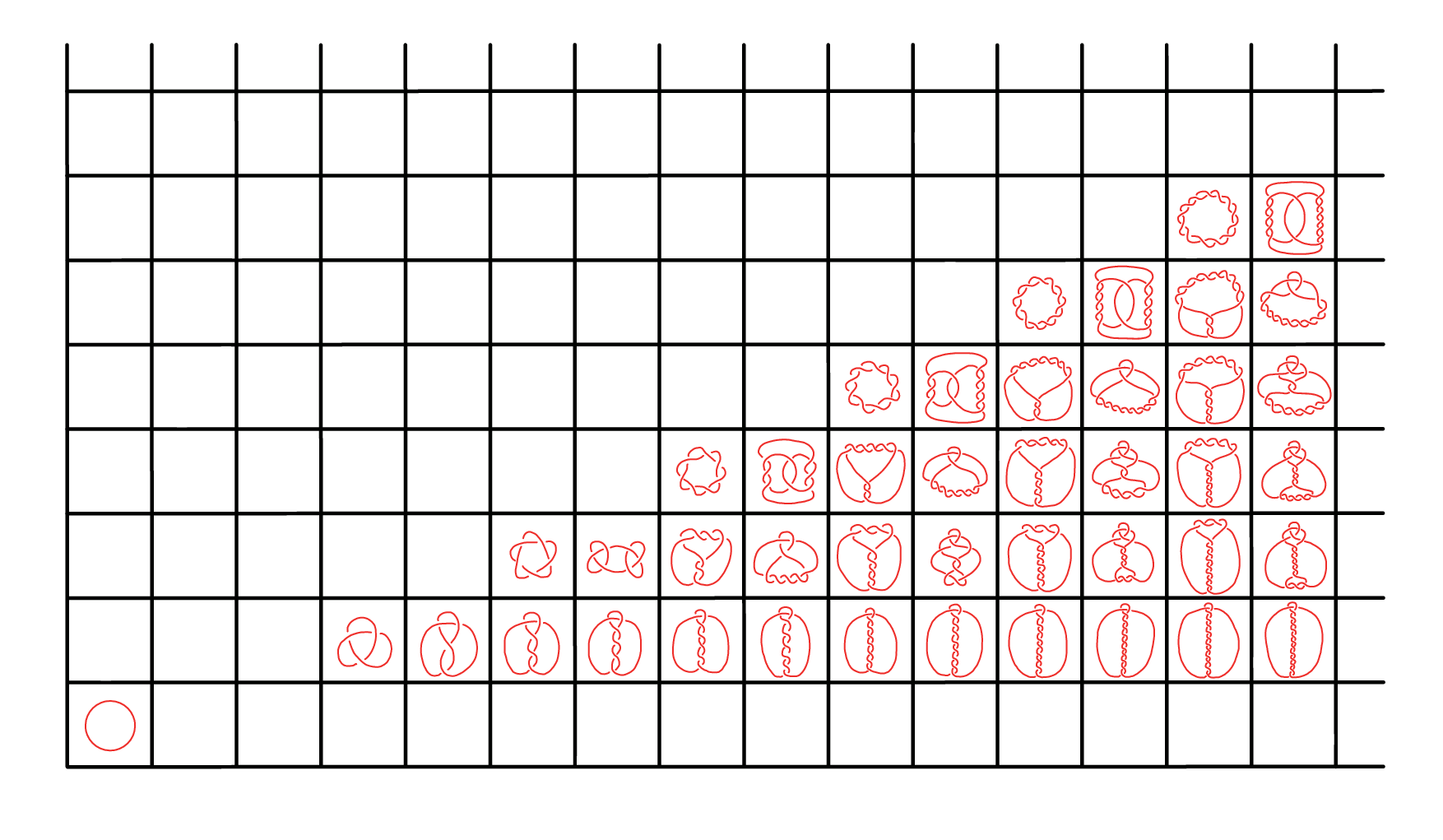}}
      \end{center}
   \caption{A choice of knots with prescribed crossing number and unknotting number}
  \label{cu-knots}
\end{figure}

\begin{figure}[htbp]
      \begin{center}
\scalebox{0.4}{\includegraphics*{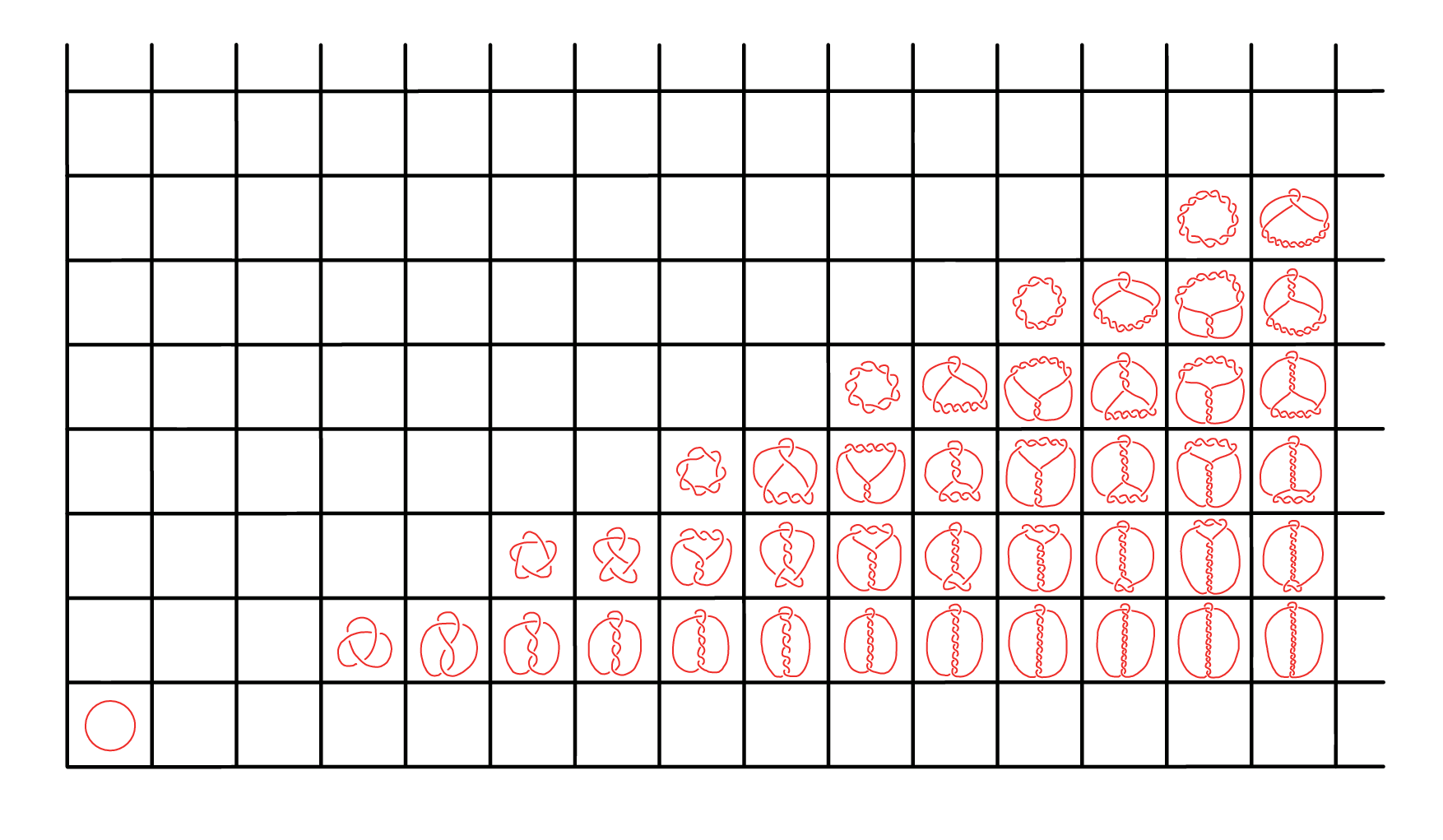}}
      \end{center}
   \caption{A choice of knots with prescribed crossing number and (canonical) genus}
  \label{cg-knots}
\end{figure}

\begin{Remark}\label{remark-cu}
{\rm
We note that the inequality $\displaystyle{u(L)\leq\frac{1}{2}(c(L))}$ holds for any link $L$. See for example \cite{Taniyama}. 
However, this inequality does not hold for planar graphs embedded in ${\mathbb S}^{3}$ \cite{A-T}. 

As mentioned above, the equality $\displaystyle{u(K)=\frac{1}{2}(c(K)-1)}$ holds if and only if $K$ is a torus knot $T(2,n)$ for some odd number $n$ \cite{Taniyama}. 
The knots $K$ with
\[
u(K)=\frac{1}{2}(c(K)-2)
\]
such as $4_{1}$, ${3_{1}}^*\#{3_{1}}^*$, $3_{1}\#{3_{1}}^{*}$ and $P(2p+1,-2,2q+1)$ are described in \cite{A-H-H}. All of them are closed $3$-braid knots. 
They contains famous knots $P(3,-2,3)$, $P(3,-2,5)$ and $P(3,-2,7)$. The first two knots are known to be all of the almost alternating torus knots \cite{Abe}. The last one, usually denoted by $P(-2,3,7)$, is known to be a key example in the theory of Dehn surgery \cite{H-O}\cite{Zhang}. 

These studies of knots with $\displaystyle{u(K)=\frac{1}{2}(c(K)-1)}$ or $\displaystyle{u(K)=\frac{1}{2}(c(K)-2)}$ are unified botany problems in knot theory. Some botany problems in knot theory are already studied. For example, genus one unknotting number one knots are doubled knots \cite{Kobayashi} \cite{S-T}, genus one three-bridge knots are pretzel knots \cite{F-O-T}, genus one tunnel number one knots are determined in \cite{Scharlemann}, and unknotting number one two-bridge knots are determined in \cite{K-M}.

We recall that the relation between $c$ and $u$ is the image of the pair $(c,u):{\mathcal K}\to({\mathbb Z}_{\geq0})^{2}$. We choose a subset ${\mathcal S}$ of ${\mathcal K}$ such that $(c,u)$ maps ${\mathcal S}$ injectively onto $(c,u)({\mathcal K})$. A choice of such ${\mathcal S}$ is illustrated in Figure \ref{cu-knots}. 

A choice for the pair $(c,g)$, that is also a choice for the pair $(c,g_{c})$, is illustrated in Figure \ref{cg-knots}. 
We note \cite{R-D} as a related study. 

}
\end{Remark}


Similar to unknotting number, the following theorem holds.


\begin{Theorem}\label{theorem-crossing-number-genus-equality}

Let $K$ be an oriented knot in ${\mathbb S}^{3}$. Then the following conditions are mutually equivalent.

\begin{enumerate}

\item[{\rm (1)}]

$\displaystyle{
g(K)=\frac{1}{2}(c(K)-1).
}$

\item[{\rm (2)}]

$\displaystyle{
g_{c}(K)=\frac{1}{2}(c(K)-1).
}$

\item[{\rm (3)}]

$K$ is a torus knot $T(2,n)$ for some odd number $n$ with $n\neq\pm1$. 

\end{enumerate}

\end{Theorem}


Ohyama showed in \cite{Ohyama} that every knot $K$ satisfies the inequality
\[
({\rm braid}-1)(K)\leq\frac{1}{2}c(K).
\]
Then the relation between crossing number and braid index minus one is described in the following theorem.

\begin{Theorem}\label{theorem-crossing-number-braid-index}
\[
(c,{\rm braid}-1)({\mathcal K})=\{(0,0)\}\cup\{(2n+1,1)\mid n\in{\mathbb Z}_{>0}\}\cup\{(x,y)\in({\mathbb Z}_{\geq2})^{2}\mid y\leq\frac{1}{2}x\}.
\]
Moreover
\[
(c,{\rm braid}-1)({\mathcal K})=(c,{\rm braid}-1)({\mathcal R}_{0}).
\]

\end{Theorem}

\vskip 5mm

\begin{figure}[htbp]
      \begin{center}
\scalebox{0.6}{\includegraphics*{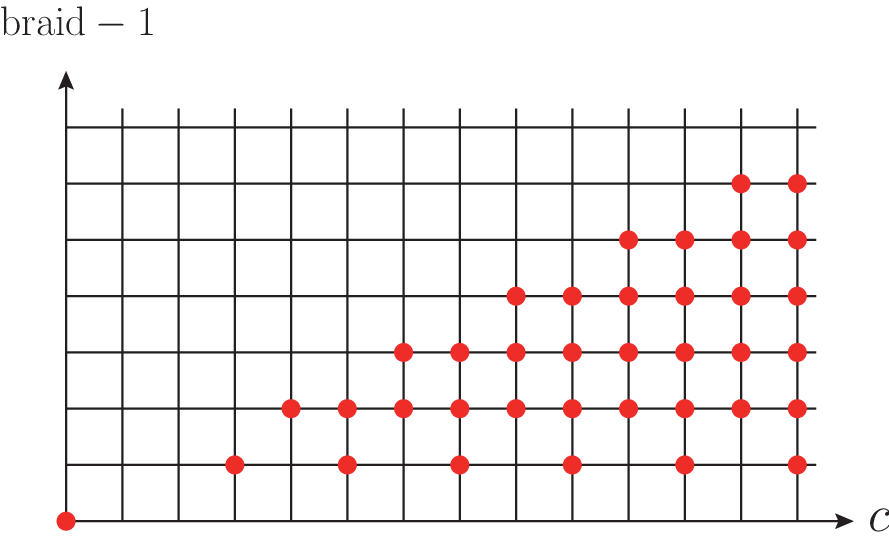}}
      \end{center}
   \caption{$(c,{\rm braid}-1)({\mathcal K})$}
  \label{cbraid-1}
\end{figure}

\begin{figure}[htbp]
      \begin{center}
\scalebox{0.4}{\includegraphics*{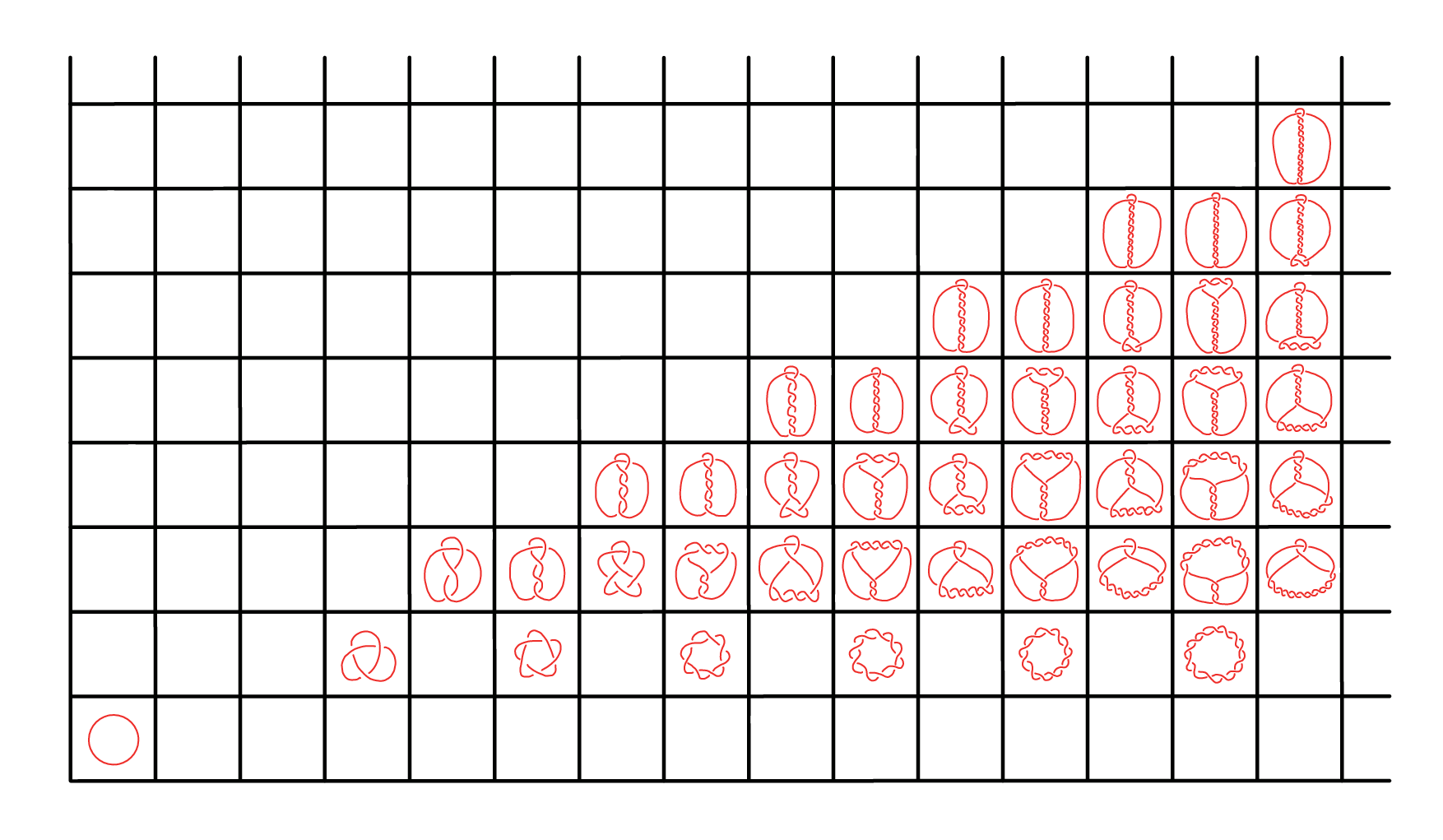}}
      \end{center}
   \caption{A choice of knots with prescribed crossing number and braid index minus one}
  \label{cbraid-1-knots}
\end{figure}

See Figure \ref{cbraid-1} and Figure \ref{cbraid-1-knots}.


\vskip 5mm

In \cite{Fox} Fox asked a question whether or not every oriented knot $K$ satisfies the inequality 
\[
{\rm bridge}(K)-1\leq\frac{1}{3}c(K).
\]
It is now called the Fox conjecture. See also \cite{Murasugi} and \cite{Ohyama}. 
Then we have the following conjecture on the relation between crossing number and bridge number minus one.

\begin{Conjecture}\label{conjecture-crossing-number-bridge-number}
\[
(c,{\rm bridge}-1)({\mathcal K})=\{(0,0)\}\cup\{(x,y)\in({\mathbb Z}_{>0})^{2}\mid y\leq\frac{1}{3}x\}.
\]

\end{Conjecture}

\vskip 5mm

\begin{figure}[htbp]
      \begin{center}
\scalebox{0.6}{\includegraphics*{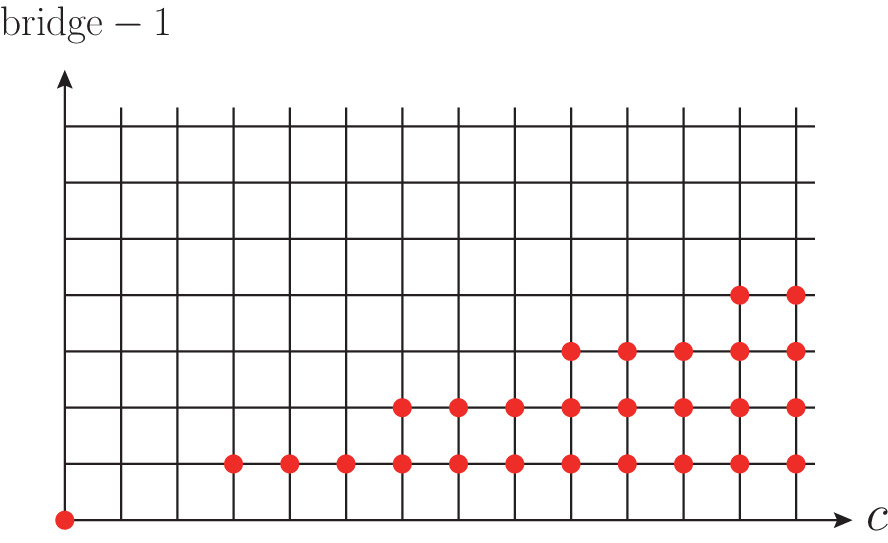}}
      \end{center}
   \caption{$(c,{\rm bridge}-1)({\mathcal K})$ (assuming Fox conjecture)}
  \label{cbridge-1}
\end{figure}

\begin{figure}[htbp]
      \begin{center}
\scalebox{0.4}{\includegraphics*{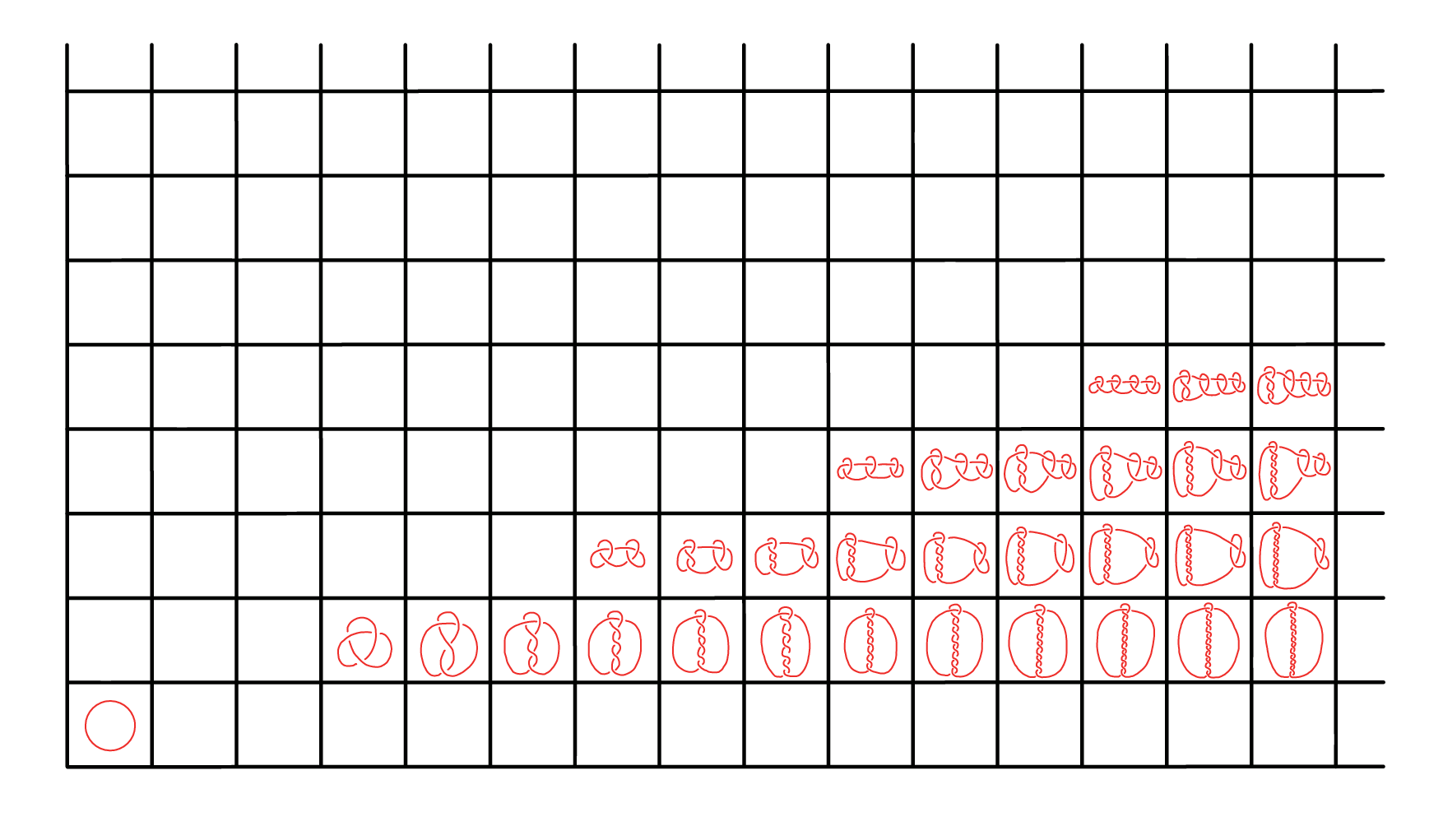}}
      \end{center}
   \caption{A choice of knots with prescribed crossing number and bridge number minus one}
  \label{cbridge-1-knots}
\end{figure}

See Figure \ref{cbridge-1} and Figure \ref{cbridge-1-knots}.


\vskip 5mm

We note that $c^{-1}(k)$ is a finite set for every $k\in{\mathbb Z}_{\geq0}$. 
Therefore every real valued function on ${\mathcal K}$ is estimated above by a function of the crossing number. 
Namely we have the following theorem.

\begin{Theorem}\label{theorem-crossing-number-estimate}

Let $f:{\mathcal K}\to{\mathbb R}$ be a function. 
Then there is a function $\varphi:{\mathbb Z}_{\geq0}\to{\mathbb R}$ such that $f(K)\leq\varphi(c(K))$ for every $K\in{\mathcal K}$.

\end{Theorem}

\vskip 5mm

\noindent{\bf Proof.} Set $\varphi(k)={\rm max}\{f(K)\mid K\in c^{-1}(k)\}$ for each $k\in{\mathbb Z}_{\geq0}$. Then we have the result. 
$\Box$

\vskip 5mm

We note here that a pair of crossing number and signature of a knot is studied in \cite{I-O-S}. 

All of the estimations $\displaystyle{u(K)\leq\frac{1}{2}(c(K)-1)}$, $\displaystyle{g(K)\leq g_{c}(K)\leq\frac{1}{2}(c(K)-1)}$,\\
$\displaystyle{({\rm braid}-1)(K)\leq\frac{1}{2}c(K)}$ and
$\displaystyle{{\rm bridge}(K)-1\leq\frac{1}{3}c(K)}$ so far are linear. It is not the case in general.

\begin{figure}[htbp]
      \begin{center}
\scalebox{0.6}{\includegraphics*{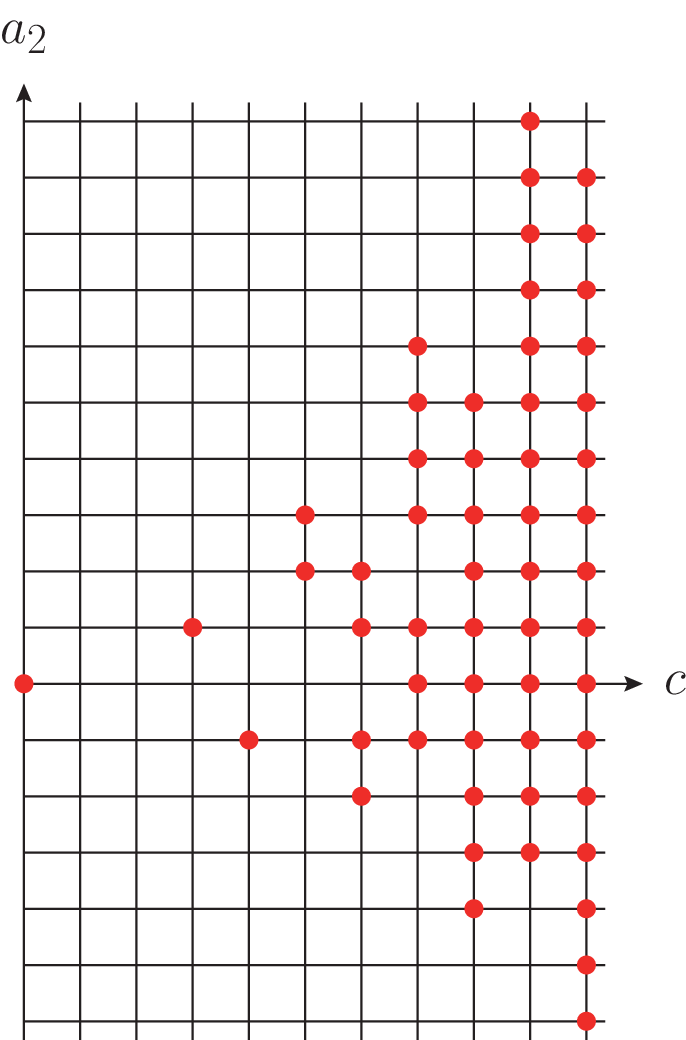}}
      \end{center}
   \caption{$(c,a_{2})({\mathcal K})$ up to $10$ crossings}
  \label{ca2}
\end{figure}

\begin{Example}\label{example-crossing-number-a2}
{\rm
Let $a_{2}:{\mathcal K}\to{\mathbb Z}$ be the second coefficient of the Conway polynomial. 
Then we have $\displaystyle{a_{2}(T(2,2k+1))=\frac{1}{2}k(k+1)}$. 
Since $c(T(2,2k+1))=2k+1$, there exists no linear estimation of $a_{2}$ by the crossing number. 
See Figure \ref{ca2}. 
}
\end{Example}

Let $R$ be a subset of ${\mathbb R}$. A map $f:{\mathcal K}\to R$ is said to be {\it crossing-order $\leq n$} if there exist real numbers $A_{0},A_{1},\cdots,A_{n}$ such that $|f(K)|\leq A_{0}+A_{1}c(K)+\cdots+A_{n}c(K)^{n}$ for every $K\in{\mathcal K}$. 
Then we see that $a_{2}$ is not crossing-order $\leq1$. 
It is shown in \cite{P-V} that $\displaystyle{|a_{2}(K)|\leq\frac{1}{8}c(K)^{2}}$ for every knot $K$. 
Therefore $a_{2}$ is crossing-order $\leq 2$. 
It is shown in \cite[Proposition 3.6]{Ito2} that $\displaystyle{a_{2}(K)\leq\frac{1}{8}(c(K)^{2}-1)}$ for every non-trivial knot $K$. 
It is also remarked in \cite[14.3]{C-D-M} that a Vassiliev invariant of order $\leq n$ has crossing-order $\leq n$. 

\begin{figure}[htbp]
      \begin{center}
\scalebox{0.5}{\includegraphics*{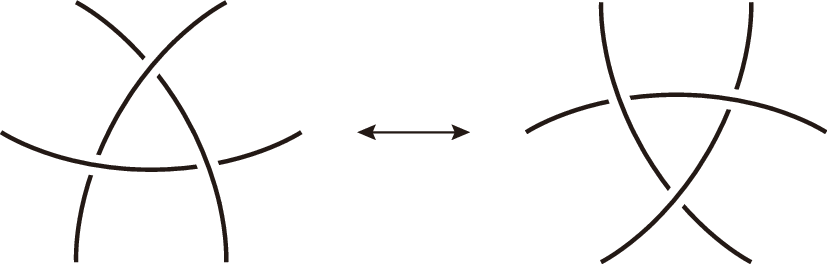}}
      \end{center}
   \caption{A delta move}
  \label{delta}
\end{figure}

A {\it delta move} is a local move as illustrated in Figure \ref{delta}. It is defined independently in \cite{Matveev} and \cite{M-N}. 
It is shown that two knots are transformed into each other by applications of delta move \cite{M-N}. 
The minimal number of applications of delta move from a given knot $K$ to a trivial knot is said to be the delta-unknotting number of $K$ and denoted by $u_{\Delta}(K)$. 
Let $u_{\Delta}:{\mathcal K}\to{\mathbb Z}_{\geq0}$ be the delta-unknotting number. It is shown in \cite{Okada} that $u_{\Delta}(K)\geq|a_{2}(K)|$ for every knot $K$. Therefore we see that $u_{\Delta}$ is not crossing-order $\leq1$. 
The next theorem shows that $u_{\Delta}$ is crossing-order $\leq2$.

\begin{Theorem}\label{theorem-crossing-number-delta-unknotting-number}

Let $K$ be an oriented knot in ${\mathbb S}^{3}$ with $c(K)\geq4$. Then 
\[
u_{\Delta}(K)\leq\frac{1}{4}(c(K)^{2}-2c(K)-3).
\]

\end{Theorem}

By definition we have ${\rm bridge}(K)\leq{\rm braid}(K)$ for every oriented knot $K$. 
Then the relation between braid index minus one and bridge number minus one is described in the following theorem.

\begin{Theorem}\label{theorem-braid-index-bridge-number}

\[
({\rm braid}-1,{\rm bridge}-1)({\mathcal K})=\{(0,0)\}\cup\{(x,y)\in({\mathbb Z}_{>0})^{2}\mid y\leq x\}.
\]

\end{Theorem}

\begin{figure}[htbp]
      \begin{center}
\scalebox{0.6}{\includegraphics*{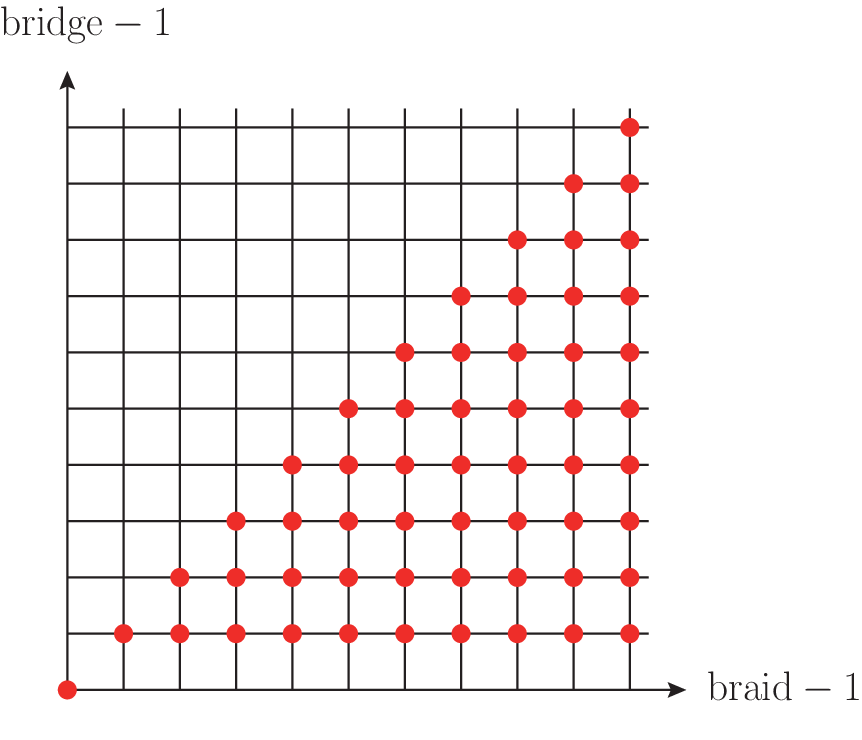}}
      \end{center}
   \caption{$({\rm braid}-1,{\rm bridge}-1)({\mathcal K})$}
  \label{braidbridge}
\end{figure}

See Figure \ref{braidbridge}.

As mentioned above $g(K)\leq g_{c}(K)$ for every oriented knot $K$. 
However the relation between canonical genus and genus is slightly different from that of braid index minus one and bridge number minus one. 
It is shown in \cite[Theorem 1.1]{Stoimenow} that no knot has $g_{c}=2$ and $g=1$. The following theorem is implicitly stated without proof in \cite{Stoimenow}.

\begin{Theorem}\label{theorem-canonical-genus-genus}

\[
(g_{c},g)({\mathcal K})=\{(0,0)\}\cup(\{(x,y)\in({\mathbb Z}_{>0})^{2}\mid y\leq x\}\setminus\{(2,1)\}).
\]

\end{Theorem}

\begin{figure}[htbp]
      \begin{center}
\scalebox{0.6}{\includegraphics*{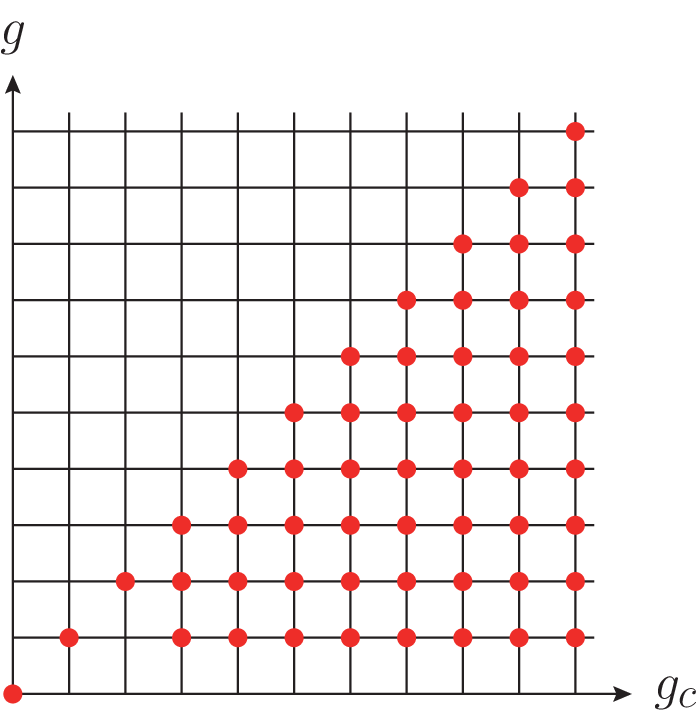}}
      \end{center}
   \caption{$(g_{c},g)({\mathcal K})$}
  \label{cgg}
\end{figure}

See Figure \ref{cgg}.

It is well-known that unknotting number and braid index are independent except $u(K)=0$ if and only if ${\rm braid}(K)=1$. 
Similar independence holds for some other pairs of oriented knot invariants. Namely we have the following theorem.

\begin{Theorem}\label{theorem-independent-pairs}

\begin{gather*}
(u,{\rm braid}-1)({\mathcal K})=(u,{\rm bridge}-1)({\mathcal K})=(g,{\rm braid}-1)({\mathcal K})\\
=(g_{c},{\rm braid}-1)({\mathcal K})=(g,{\rm bridge}-1)({\mathcal K})=\{(0,0)\}\cup({\mathbb Z}_{>0})^{2}.
\end{gather*}

\end{Theorem}

\begin{figure}[htbp]
      \begin{center}
\scalebox{0.6}{\includegraphics*{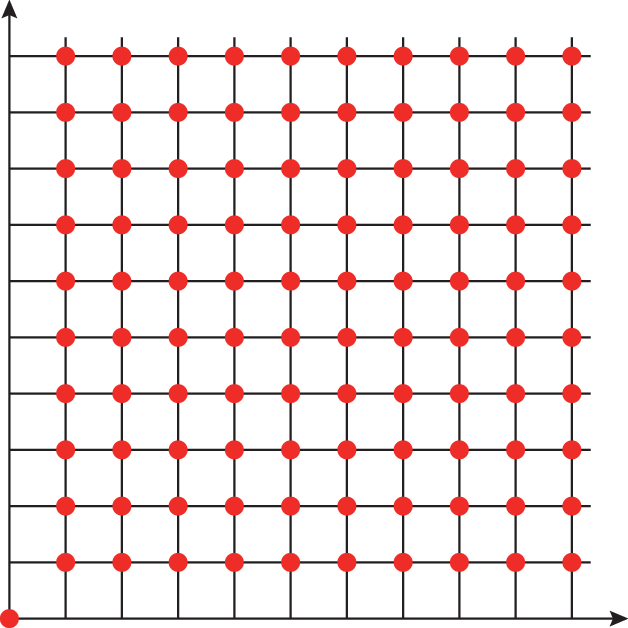}}
      \end{center}
   \caption{$(u,{\rm braid}-1)({\mathcal K})=(u,{\rm bridge}-1)({\mathcal K})
   \newline 
   =(g,{\rm braid}-1)({\mathcal K})=(g_{c},{\rm braid}-1)({\mathcal K})=(g,{\rm bridge}-1)({\mathcal K})$}
  \label{independent}
\end{figure}

\begin{figure}[htbp]
      \begin{center}
\scalebox{0.4}{\includegraphics*{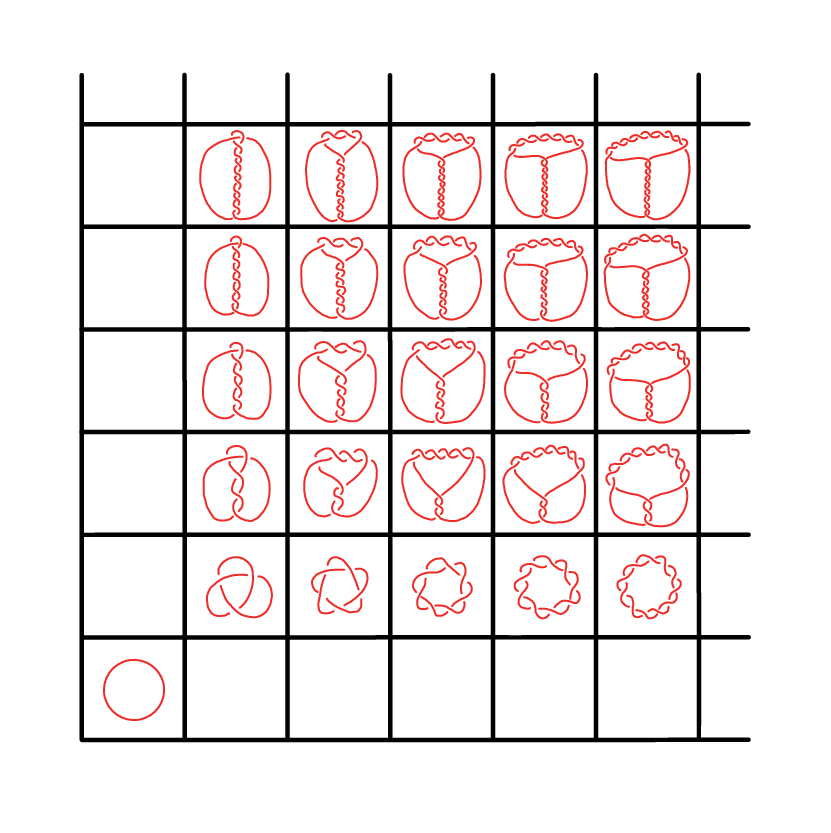}}
      \end{center}
   \caption{A choice of knots with prescribed canonical genus and braid index minus one}
  \label{gbraid-1-knots}
\end{figure}

See Figure \ref{independent} and Figure \ref{gbraid-1-knots}.

%


We mainly consider crossing number, unknotting number, braid index minus one, bridge number minus one, genus and canonical genus. 
There are $15$ unordered pairs of them. We note that $12$ of them are decided as above. 
The rest are $(u,g)$, $(u,g_{c})$ and $(g_{c},{\rm bridge}-1)$. 
It seems to us that they are a bit difficult to determine. 
To the best of our knowledge, no knot $K$ with $g(K)=1$ and $u(K)\geq4$ is known. 
It is only known that $g(P(3,3,3))=1$ and $u(P(3,3,3))=3$ \cite{Owens}. 
See also \cite{L-L}. 
It is shown in \cite[Theorem 2.1]{Stoimenow2} that a canonical genus $1$ knot is a $2$-bridge knot with Conway notation $C(2p,2q)$, or a pretzel knot $P(2p+1,2q+1,2r+1)$. Thus there exists no knot $K$ of $g_{c}(K)=1$ and ${\rm bridge}(K)\geq4$. 
In general we have the following theorem.

\begin{Theorem}\label{theorem-canonical-genus-bridge-number}

Let $k$ be a positive integer. Then there exists a positive integer $m$ such that every oriented knot $K$ with $g_{c}(K)=k$ satisfies ${\rm bridge}(K)\leq m$.

\end{Theorem}

\section{Crossing number and unknotting number, genus or canonical genus}\label{section-cuggc} 
%


Let $K$ be an oriented knot in ${\mathbb S}^{3}$ and $D$ a knot diagram of $K$. 
The number of crossings of $D$ is denoted by $c(D)$. 
A canonical Seifert surface of $K$ obtained from $D$ by Seifert's algorithm is denoted by $F(D)$. 
We denote the number of Seifert circles obtained from $D$ by $s(D)$.

\vskip 5mm

\noindent{\bf Proof of Theorem \ref{theorem-crossing-number-unknotting-number-genus-canonical-genus}.} 
We divide the proof into the following two claims. 

\vskip 3mm
\noindent
Claim 1.
\begin{gather*}
(c,u)({\mathcal K})=(c,u)(\{{3_{1}}^*\#{3_{1}}^*\})\cup(c,u)({\mathcal P}_{0})\\
=\{(0,0)\}\cup\{(x,y)\in({\mathbb Z}_{>0})^{2}\mid y\leq\frac{1}{2}(x-1)\}.
\end{gather*}

\vskip 3mm
\noindent
Claim 2.
\begin{gather*}
(c,g)({\mathcal K})=(c,g)({\mathcal R}_{0})=(c,g_{c})({\mathcal K})=(c,g_{c})({\mathcal R}_{0})\\
=\{(0,0)\}\cup\{(x,y)\in({\mathbb Z}_{>0})^{2}\mid y\leq\frac{1}{2}(x-1)\}.
\end{gather*}

First we will show Claim 1. For a knot $K$, three conditions that $K$ is trivial, $c(K)=0$ and $u(K)=0$ are mutually equivalent. Therefore $K$ is trivial if and only if $(c,u)(K)=(0,0)$, and if $K$ is non-trivial then $(c,u)(K)\in({\mathbb Z}_{>0})^{2}$. 
As mentioned above, every nontrivial knot $K$ satisfies $\displaystyle{u(K)\leq\frac{1}{2}(c(K)-1)}$. 
Therefore we have 
\[
\{(0,0)\}\cup\{(x,y)\in({\mathbb Z}_{>0})^{2}\mid y\leq\frac{1}{2}(x-1)\}\supset(c,u)({\mathcal K}).
\]
By definition we have
\[
(c,u)({\mathcal K})\supset(c,u)(\{{3_{1}}^*\#{3_{1}}^*\})\cup(c,u)({\mathcal P}_{0}).
\]
Therefore it is sufficient to show
\[
(c,u)(\{{3_{1}}^*\#{3_{1}}^*\})\cup(c,u)({\mathcal P}_{0})\supset\{(0,0)\}\cup\{(x,y)\in({\mathbb Z}_{>0})^{2}\mid y\leq\frac{1}{2}(x-1)\}.
\]
As mentioned above, $(c,u)(0_{1})=(0,0)$. Therefore $(0,0)\in(c,u)({\mathcal P}_{0})$. 
We set
\[
{\mathcal O}=\{(x,y)\in({\mathbb Z}_{>0})^{2}\mid y\leq\frac{1}{2}(x-1),\ x\equiv 1\ (\bmod\ 2)\}
\]
and
\[
{\mathcal E}=\{(x,y)\in({\mathbb Z}_{>0})^{2}\mid y\leq\frac{1}{2}(x-1),\ x\equiv 0\ (\bmod\ 2)\}.
\]
Then we have 
\[
{\mathcal O}\cup{\mathcal E}=\{(x,y)\in({\mathbb Z}_{>0})^{2}\mid y\leq\frac{1}{2}(x-1)\}.
\] 
Therefore it is sufficient to show ${\mathcal O}\subset(c,u)(\{{3_{1}}^*\#{3_{1}}^*\})\cup(c,u)({\mathcal P}_{0})$ and\\ 
${\mathcal E}\subset(c,u)(\{{3_{1}}^*\#{3_{1}}^*\})\cup(c,u)({\mathcal P}_{0})$. 

Each element of ${\mathcal O}$ can be expressed as $(2k+1,i)$ where $k$ and $i$ are positive integers with $i\leq k$. 
We will show $(c,u)(C(2k-2i+1,2i))=(2k+1,i)$. 
We note that $C(2k-2i+1,2i)$ is a $2$-bridge knot. In particular, when $i=1$, $C(2k-1,2)$ is a $(2k-1)$-twist knot and 
when $i=k$, $C(1,2k)=C(2k+1)$ is a torus knot $T(2,-(2k+1))$. 
Since $C(2k-2i+1,2i)$ is alternating and both $2k-2i+1$ and $2i$ are positive, we have $c(C(2k-2i+1,2i))=2k-2i+1+2i=2k+1$. 
Let $\sigma(K)$ be the signature of a knot $K$. 
It is well-known that $\sigma(C(1,2k))=\sigma(T(2,-(2k+1)))=2k$ and $u(C(1,2k))=u(T(2,-(2k+1)))=k$. 
We see that $C(1,2i)=T(2,-(2i+1))$ is obtained from $C(2k-2i+1,2i)$ by changing $k-i$ negative crossings to positive crossings. 
Then by the inequality 
\[
\sigma
\left(
\begin{minipage}{14pt}
\scalebox{0.153}{\includegraphics{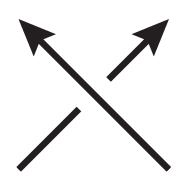}}
\end{minipage}
\right)
\leq
\sigma
\left(
\begin{minipage}{14pt}
\scalebox{0.153}{\includegraphics{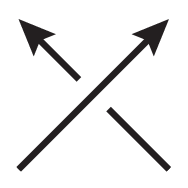}}
\end{minipage}
\right)
\]
observed in \cite{Giller} we have $\sigma(C(2k-2i+1,2i))\geq\sigma(C(1,2i))=2i$. 
Therefore we have $u(C(2k-2i+1,2i))\geq i$. 
By changing $i$ crossings of $C(2k-2i+1,2i)$ we have a trivial knot $C(2k-2i+1,0)$. 
Therefore $u(C(2k-2i+1,2i))\leq i$ and we have $u(C(2k-2i+1,2i))=i$. 
Thus we have $(c,u)(C(2k-2i+1,2i))=(2k+1,i)$ as intended. 
Therefore $(2k+1,i)$ is an element of $(c,u)({\mathcal R})\subset(c,u)({\mathcal P}_{0})$. 
Thus we have ${\mathcal O}\subset(c,u)(\{{3_{1}}^*\#{3_{1}}^*\})\cup(c,u)({\mathcal P}_{0})$. 

Each element of ${\mathcal E}$ can be expressed as $(2k,i)$ where $k$ and $i$ are positive integers with $i<k$. 
First we consider the case $i=k-1$. 
For $i=1$, we have\\ 
$(2k,i)=(4,1)=(c,u)(4_{1})$ where $4_{1}=C(2,2)$ is the figure-eight knot. 
For $i=2$, we have $(2k,i)=(6,2)$. All $6$-crossing prime knots have unknotting number $1$. Namely, $u(6_1)=u({6_1}^{*})=u(6_2)=u({6_2}^{*})=u(6_3)=1$. 
Then we see that $(c,u)(K)=(6,2)$ if and only if $K$ is ${3_{1}}^*\#{3_{1}}^*$, $3_{1}\#{3_{1}}^{*}$ or $3_{1}\#3_{1}$. 
Thus $(6,2)$ is not an element of $(c,u)({\mathcal P}_{0})$ but an element of $(c,u)(\{{3_{1}}^*\#{3_{1}}^*\})$. 

For $i\geq3$, we consider a pretzel knot $P(2p+1,-2,2q+1)$ where $p$ and $q$ are positive integers with $p+q=i-1$. 
We will show
\[
(c,u)(P(2p+1,-2,2q+1))=(2i+2,i)=(2k,i).
\]
We note that the standard diagram of $P(2p+1,-2,2q+1)$ is a diagram with $2p+2q+4=2i+2$ crossings. 
Therefore $c(P(2p+1,-2,2q+1))\leq2i+2$. Then 
\[
u(P(2p+1,-2,2q+1))\leq\frac{1}{2}(c(P(2p+1,-2,2q+1))-1)\leq i+\frac{1}{2}.
\]
Therefore $u(P(2p+1,-2,2q+1))\leq i$. 
Let $S(K)$ be the Rasmussen invariant of an oriented knot $K$ defined in \cite{Rasmussen}. 
Suppose that $D$ is a positive diagram of $K$. Then it is shown in \cite{Rasmussen} that $S(K)=c(D)-s(D)+1$. 
Since the standard diagram of $P(2p+1,-2,2q+1)$ is a positive $3$-braid diagram, we have $S(P(2p+1,-2,2q+1))=(2i+2)-3+1=2i$. 
Then by \cite{Rasmussen} we have $u(P(2p+1,-2,2q+1))\geq i$. 
Thus we have 
$u(P(2p+1,-2,2q+1))=i$. 
Therefore we have $c(P(2p+1,-2,2q+1))\geq2i+1$. 
Suppose $c(P(2p+1,-2,2q+1))=2i+1$. Then by \cite[Theorem 1.4]{Taniyama} we see that $P(2p+1,-2,2q+1)$ is a torus knot $T(2,2i+1)$. It is shown in \cite{Shinohara} that 
\[
\sigma(P(2p+1,-2,2q+1))=
\left\{ \begin{array}{ll}
-2i & (i<4)\\
-2i+2 & (i\geq4).\\
\end{array} \right.
\]
It is shown in \cite{Rasmussen} that an alternating knot $K$ satisfies $S(K)=-\sigma(K)$. 
Then we see that $P(2p+1,-2,2q+1)$ with $i\geq4$ is non-alternating. This contradicts that $T(2,2i+1)$ is alternating. The knots $P(2p+1,-2,2q+1)$ with $i<4$ are $P(3,-2,3)$. This is a torus knot $T(3,4)$. Thus we have $c(P(2p+1,-2,2q+1))=2i+2$. 
In summary, we have $(c,u)(P(2p+1,-2,2q+1))=(2i+2,i)$ as intended. 
Suppose that $P(2p+1,-2,2q+1)$ is a composite knot $J\#K$. Since $({\rm braid}-1)(P(2p+1,-2,2q+1))=2$, we have $({\rm braid}-1)(J)=({\rm braid}-1)(K)=1$. Then we have $J=T(2,a)$ and $K=T(2,b)$ for some odd integers $a$ and $b$. Then $J\#K$ is an alternating knot. This contradicts the fact that $P(2p+1,-2,2q+1)$ is non-alternating. 
Thus we see that $P(2p+1,-2,2q+1)$ is a prime knot. 

Next we consider the case $i\leq k-2$. 
We have $(c,u)(C(2k-2,2))=(2k,1)$ where $C(2k-2,2)$ is a $(2k-2)$-twist knot. 
Therefore we consider the case $2\leq i\leq k-2$. Then we have $k\geq4$. 
We will show that $(c,u)(C(2k-3,1,2))=(2k,k-2)$. 
We will also show for $2\leq i\leq k-3$ that $(c,u)(C(2i-1,2k-2i-4,2,1,2))=(2k,i)$. 
Since $2k-3$, $1$ and $2$ are positive we have $c(C(2k-3,1,2))=2k-3+1+2=2k$. 
By \cite{Shinohara} we have $\sigma(C(2k-3,1,2))=2k-4$. 
Then we have $u(C(2k-3,1,2))\geq k-2$. 
By changing $k-2$ crossings of $C(2k-3,1,2)$ we have a trivial knot $C(3,-1,2)$. 
Thus we have $u(C(2k-3,1,2))=k-2$. Thus we have $(c,u)(C(2k-3,1,2))=(2k,k-2)$. 
Since $2\leq i\leq k-3$, all of $2i-1$, $2k-2i-4$, $2$, $1$ and $2$ are positive. 
Therefore we have $c(C(2i-1,2k-2i-4,2,1,2))=2i-1+2k-2i-4+2+1+2=2k$. 
By changing $k-i-2$ negative crossings of $C(2i-1,2k-2i-4,2,1,2)$ to positive crossings, we have an oriented knot $C(2i-1,0,2,1,2)=C(2i+1,1,2)$. Then by the inequality stated above we have $\sigma(C(2i-1,2k-2i-4,2,1,2))\geq\sigma(C(2i+1,1,2))=2i$. 
Thus we have $u(C(2i-1,2k-2i-4,2,1,2))\geq i$. 
By changing $i$ crossings of $C(2i-1,2k-2i-4,2,1,2)$, we have a trivial knot $C(1,2k-2i-4,2,-1,2)$. 
Therefore $u(C(2i-1,2k-2i-4,2,1,2))\leq i$ and we have 
$(c,u)(C(2i-1,2k-2i-4,2,1,2))=(2k,i)$.
Thus we have 
\[
(c,u)(\{{3_{1}}^*\#{3_{1}}^*\})\cup(c,u)({\mathcal P}_{0})\supset\{(0,0)\}\cup\{(x,y)\in({\mathbb Z}_{>0})^{2}\mid y\leq\frac{1}{2}(x-1)\}
\]
and we have shown Claim 1. 

Next we will show Claim 2. 
It is sufficient to show 
\[
(c,g)({\mathcal R}_{0})\subset(c,g)({\mathcal K})\subset\{(0,0)\}\cup\{(x,y)\in({\mathbb Z}_{>0})^{2}\mid y\leq\frac{1}{2}(x-1)\}\subset(c,g)({\mathcal R}_{0})
\]
and
\[
(c,g_{c})({\mathcal R}_{0})\subset(c,g_{c})({\mathcal K})\subset\{(0,0)\}\cup\{(x,y)\in({\mathbb Z}_{>0})^{2}\mid y\leq\frac{1}{2}(x-1)\}\subset(c,g_{c})({\mathcal R}_{0}).
\]
Since ${\mathcal R}_{0}\subset{\mathcal K}$ we have $(c,g)({\mathcal R}_{0})\subset(c,g)({\mathcal K})$ and $(c,g_{c})({\mathcal R}_{0})\subset(c,g_{c})({\mathcal K})$. 

We will show
\[
(c,g)({\mathcal K})\subset\{(0,0)\}\cup\{(x,y)\in({\mathbb Z}_{>0})^{2}\mid y\leq\frac{1}{2}(x-1)\}
\]
and
\[
(c,g_{c})({\mathcal K})\subset\{(0,0)\}\cup\{(x,y)\in({\mathbb Z}_{>0})^{2}\mid y\leq\frac{1}{2}(x-1)\}.
\]
For a knot $K$, four conditions that $K$ is trivial, $c(K)=0$, $g(K)=0$ and $g_{c}(K)=0$ are mutually equivalent. Therefore if $K$ is trivial, then $(c,g)(K)=(c,g_{c})(K)=(0,0)$, and if $K$ is non-trivial, then both $(c,g)(K)$ and $(c,g_{c})(K)$ are elements of $({\mathbb Z}_{>0})^{2}$. 
Let $K$ be a non-trivial knot. By definition we have $g(K)\leq g_{c}(K)$. Let $D$ be a diagram of $K$ with $c(D)=c(K)$. 
Then we have $\chi(F(D))=s(D)-c(D)$. Since $\chi(F(D))=1-2g(F(D))$ and $s(D)\geq2$, we have 
\[
g(F(D))=\frac{1}{2}(1-s(D)+c(D))\leq\frac{1}{2}(1-2+c(D))=\frac{1}{2}(c(D)-1).
\]
Since $g_{c}(K)\leq g(F(D))$ and $c(D)=c(K)$ we have 
\[
g(K)\leq g_{c}(K)\leq\frac{1}{2}(c(K)-1).
\]
Therefore both $(c,g)(K)=(c(K),g(K))$ and $(c,g_{c})(K)=(c(K),g_{c}(K))$ are elements of $\displaystyle{\{(x,y)\in({\mathbb Z}_{>0})^{2}\mid y\leq\frac{1}{2}(x-1)\}}$. Thus we have shown
\[
(c,g)({\mathcal K})\subset\{(0,0)\}\cup\{(x,y)\in({\mathbb Z}_{>0})^{2}\mid y\leq\frac{1}{2}(x-1)\}
\]
and
\[
(c,g_{c})({\mathcal K})\subset\{(0,0)\}\cup\{(x,y)\in({\mathbb Z}_{>0})^{2}\mid y\leq\frac{1}{2}(x-1)\}.
\]
Finally we will show
\[
\{(0,0)\}\cup\{(x,y)\in({\mathbb Z}_{>0})^{2}\mid y\leq\frac{1}{2}(x-1)\}\subset(c,g)({\mathcal R}_{0})
\]
and
\[
\{(0,0)\}\cup\{(x,y)\in({\mathbb Z}_{>0})^{2}\mid y\leq\frac{1}{2}(x-1)\}\subset(c,g_{c})({\mathcal R}_{0}).
\]
Since $(c,g)(0_{1})=(c,g_{c})(0_{1})=(0,0)$ we see that $(0,0)$ is an element of both $(c,g)({\mathcal R}_{0})$ and $(c,g_{c})({\mathcal R}_{0})$. 
Let $(x,y)$ be an element of ${\mathcal O}$. 
Let $K$ be a $2$-bridge knot with Conway notation $C(2y,x-2y)$. 
Let $D$ be a knot diagram corresponding to this Conway notation. 
Since $D$ is a reduced alternating diagram of $K$ with $c(D)=2y+(x-2y)=x$, we have $c(K)=x$. 
We see $s(D)=x-2y+1$, $\chi(F(D))=-2y+1$ and $g(F(D))=y$. 
It is known that the genus of an alternating knot is equal to the genus of a canonical Seifert surface obtained from an alternating diagram. See \cite{Murasugi3} \cite{Crowell} \cite{Gabai}. 
Then we have $g(K)=g_{c}(K)=g(F(D))=y$. 
Therefore we have $(c,g)(K)=(c,g_{c})(K)=(x,y)$. 
Thus $(x,y)$ is an element of both $(c,g)({\mathcal R}_{0})$ and $(c,g_{c})({\mathcal R}_{0})$. 
Let $(x,y)$ be an element of ${\mathcal E}$. 
Let $K$ be a $2$-bridge knot with Conway notation $C(-2,-x+2y+1,-2y+1)$. 
Let $D$ be a knot diagram corresponding to this Conway notation. 
Since $D$ is a reduced alternating diagram of $K$ with $c(D)=-(-2+(-x+2y+1)+(-2y+1))=x$, we have $c(K)=x$. 
We see $s(D)=x-2y+1$, $\chi(F(D))=-2y+1$ and $g(F(D))=y$. 
Then we have $g(K)=g_{c}(K)=g(F(D))=y$. 
Therefore we have $(c,g)(K)=(c,g_{c})(K)=(x,y)$. 
Thus $(x,y)$ is an element of both $(c,g)({\mathcal R}_{0})$ and $(c,g_{c})({\mathcal R}_{0})$. 
This completes the proof. 
$\Box$

\vskip 5mm

\noindent{\bf Proof of Theorem \ref{theorem-crossing-number-genus-equality}.} 
Suppose 
$\displaystyle{
g(K)=\frac{1}{2}(c(K)-1).
}$
Since $g(K)\leq g_{c}(K)$ and 
$\displaystyle{
g_{c}(K)\leq\frac{1}{2}(c(K)-1)
}$
by Theorem \ref{theorem-crossing-number-unknotting-number-genus-canonical-genus}, we have
$\displaystyle{
g_{c}(K)=\frac{1}{2}(c(K)-1).
}$
Thus (1) implies (2). 
Suppose 
$\displaystyle{
g_{c}(K)=\frac{1}{2}(c(K)-1).
}$
Let $D$ be a diagram of $K$ with $c(D)=c(K)$. Then $g(F(D))\geq g_{c}(K)$. 
Therefore $\chi(F(D))=1-2g(F(D))\leq1-2g_{c}(K)=2-c(K)$. Since $\chi(F(D))=s(D)-c(D)=s(D)-c(K)$ we have $s(D)\leq2$. 
Since $K$ is not a trivial knot, $K$ is a torus knot $T(2,n)$ for some odd number $n$ with $n\neq\pm1$. 
Thus (2) implies (3). 
Suppose that $K$ is a torus knot $T(2,n)$ for some odd number $n$ with $n\neq\pm1$. 
Then we have 
$\displaystyle{
g(K)=g_{c}(K)=\frac{1}{2}(|n|-1)=\frac{1}{2}(c(K)-1).
}$
Thus (3) implies (1). 
$\Box$

\section{Crossing number and braid index minus one}\label{section-cbraid} 

\noindent{\bf Proof of Theorem \ref{theorem-crossing-number-braid-index}.} 
We set
\[
S=\{(0,0)\}\cup\{(2n+1,1)\mid n\in{\mathbb Z}_{>0}\}\cup\{(x,y)\in({\mathbb Z}_{\geq2})^{2}\mid y\leq\frac{1}{2}x\}.
\]
First we will show
\[
S\supset(c,{\rm braid}-1)({\mathcal K}).
\]
For a knot $K$, three conditions that $K$ is trivial, $c(K)=0$ and $({\rm braid}-1)(K)=0$ are mutually equivalent. Therefore $K$ is trivial if and only if $(c,{\rm braid}-1)(K)=(0,0)$, and if $K$ is non-trivial, then $(c,{\rm braid}-1)(K)\in({\mathbb Z}_{>0})^{2}$. 
Suppose that $({\rm braid}-1)(K)=1$. Then $K$ is a $(2,p)$-torus knot for some odd number $p$ with $|p|\geq3$. 
Therefore $(c,{\rm braid}-1)(K)=(|p|,1)$. 
As mentioned above every nontrivial knot $K$ satisfies $\displaystyle{({\rm braid}-1)(K)\leq\frac{1}{2}c(K)}$ \cite{Ohyama}. 
Therefore we have 
\[
S\supset(c,{\rm braid}-1)({\mathcal K}).
\]
Next we will show
\[
(c,{\rm braid}-1)({\mathcal R}_{0})\supset S.
\]
Since $(c,{\rm braid}-1)(0_{1})=(0,0)$, we have $(0,0)\in(c,{\rm braid}-1)({\mathcal R}_{0})$. 
The braid index of a $2$-bridge knot is determined in \cite[Theorem B]{Murasugi2}. 
We note that negative continued fractions are used in \cite[Theorem B]{Murasugi2}. 
Then we see
\[
(c,{\rm braid}-1)(C(2a,2b-1))=(2a+2b-1,b)
\]
and
\[
(c,{\rm braid}-1)(C(-2,-2a+1,-2b+1))=(2a+2b,a+1)
\]
for positive integers $a$ and $b$.
Let ${\mathcal S}$ be a subset of ${\mathcal R}_{0}$ defined by
\[
{\mathcal S}=\{0_{1}\}\cup\{C(2a,2b-1)\mid a,b\in{\mathbb Z}_{>0}\}\cup\{C(-2,-2a+1,-2b+1)\mid a,b\in{\mathbb Z}_{>0}\}.
\]
Then we see that $(c,{\rm braid}-1)$ maps ${\mathcal S}$ injectively onto $S$. See Figure \ref{cbraid-1-knots}.
Thus we have
\[
(c,{\rm braid}-1)({\mathcal R}_{0})\supset S.
\]
Since
\[
(c,{\rm braid}-1)({\mathcal K})\supset(c,{\rm braid}-1)({\mathcal R}_{0})
\]
we have
\[
S\supset(c,{\rm braid}-1)({\mathcal K})\supset(c,{\rm braid}-1)({\mathcal R}_{0})\supset S.
\]
Therefore we have
\[
S=(c,{\rm braid}-1)({\mathcal K})=(c,{\rm braid}-1)({\mathcal R}_{0}).
\]
$\Box$

\section{Crossing number and bridge number minus one}\label{section-cbridge} 
Here we show the following easy half of Conjecture \ref{conjecture-crossing-number-bridge-number}.

\begin{Proposition}\label{proposition-crossing-number-bridge-number}
\[
(c,{\rm bridge}-1)({\mathcal K})\supset\{(0,0)\}\cup\{(x,y)\in({\mathbb Z}_{>0})^{2}\mid y\leq\frac{1}{3}x\}.
\]

\end{Proposition}

\vskip 5mm

\noindent{\bf Proof.} 
We note that $(c,{\rm bridge}-1)(0_{1})=(0,0)$. 
Let $x$ and $y$ be positive integers with $\displaystyle{y\leq\frac{1}{3}x}$. 
Let $C(2,x-3y+1)$ be a twist knot. Then we have
\[
(c,{\rm bridge}-1)(C(2,x-3y+1)\#(y-1)\cdot{3_{1}}^*)=(x,y).
\]
This completes the proof. 
$\Box$

\vskip 5mm

See Figure \ref{cbridge-1-knots}.

\vskip 5mm

\section{Crossing number and delta-unknotting number}\label{section-cdelta} 

Let $D$ be a diagram of a knot. 
A crossing $x$ of $D$ is said to be {\it outermost} if at least one of two closed curves obtained from $D$ by smoothing $x$ is simple. 

\begin{Lemma}\label{lemma-crossing-change-by-delta-moves} Let $D$ be a diagram of a knot $K$ and $x$ a crossing of $D$. A crossing change of $K$ at $x$ is realized by at most $\displaystyle{\left\lfloor\frac{c(D)+1}{2}\right\rfloor}$-times applications of delta move on $K$. 
If $x$ is outermost, then it is realized by at most $\displaystyle{\left\lfloor\frac{c(D)-1}{2}\right\rfloor}$-times applications of delta move on $K$. 
\end{Lemma}

\begin{figure}[htbp]
      \begin{center}
\scalebox{0.6}{\includegraphics*{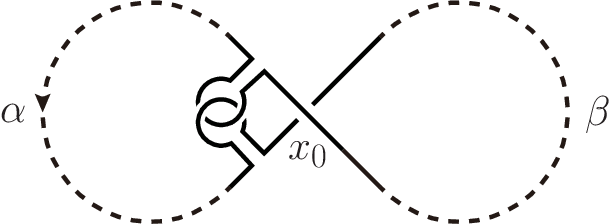}}
      \end{center}
   \caption{$K$ is a band sum of a Hopf link and $K'$}
  \label{Hopf-band}
\end{figure}

\begin{figure}[htbp]
      \begin{center}
\scalebox{0.6}{\includegraphics*{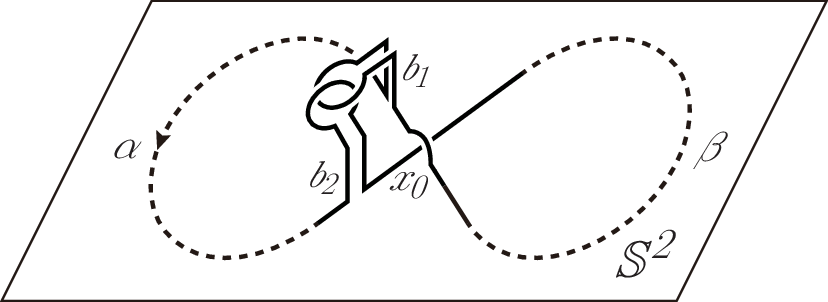}}
      \end{center}
   \caption{A Hopf link and bands form an overbridge}
  \label{Hopf-band2}
\end{figure}

\begin{figure}[htbp]
      \begin{center}
\scalebox{0.6}{\includegraphics*{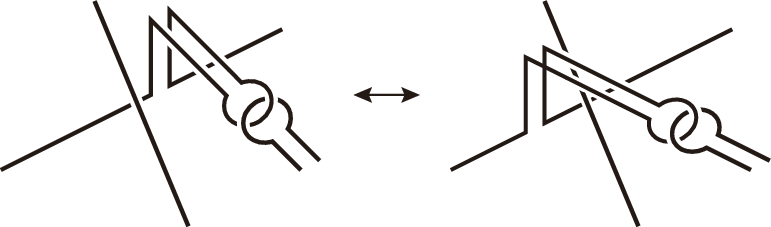}}
      \end{center}
   \caption{A clasp-leaps-over-a-hurdle move}
  \label{hurdle}
\end{figure}

\noindent{\bf Proof.} We show the case that $x$ is a positive crossing. The mirror image move of a delta move is again a delta move \cite{M-N}. Therefore the case $x$ is a negative crossing is shown similarly. 
Let $D'$ be a knot diagram obtained from $D$ by changing over/under crossing information at $x$. 
We denote the changed crossing by $x'$. Namely $x'$ is a negative crossing of $D'$. Let $K'$ be a knot represented by $D'$. 
We will show that $K'$ is obtained from $K$ by at most $\displaystyle{\left\lfloor\frac{c(D)+1}{2}\right\rfloor}$-times applications of delta move, and in case $x$ is outermost, at most $\displaystyle{\left\lfloor\frac{c(D)-1}{2}\right\rfloor}$-times applications. 
We set $x_{0}=x'$. 
Let ${\mathcal C}(D')=\{x_{0},x_{1},\cdots,x_{n}\}$ be the set of all crossings of $D'$ where $n=c(D')-1=c(D)-1$. 
Let $x_{i}^{+}$ and $x_{i}^{-}$ be the preimage points of $x_{i}$ on $K'$ such that $x_{i}^{+}$ is over $x_{i}^{-}$ on $D'$ for each $i\in\{0,1,\cdots,n\}$. 
The points $x_{0}^{+},x_{1}^{+},\cdots,x_{n}^{+}$ (resp. $x_{0}^{-},x_{1}^{-},\cdots,x_{n}^{-}$) are said to be {\it over-crossings} (resp. {\it under-crossings}). 
The point $x_{i}^{+}$ (resp. $x_{i}^{-}$) is said to be a {\it partner} of $x_{i}^{-}$ (resp. $x_{i}^{+}$) and we denote it by $x_{i}^{+}={\rm partner}(x_{i}^{-})$ (resp. $x_{i}^{-}={\rm partner}(x_{i}^{+})$) for $i\in\{0,1,\cdots,n\}$. 
Let $\alpha$ and $\beta$ be simple arcs of $K'$ such that $K'=\alpha\cup\beta$ and $\alpha\cap\beta=\partial\alpha=\partial\beta=\{x_{0}^{+},x_{0}^{-}\}$. 
We may suppose without loss of generality that $|\alpha\cap\{x_{1}^{+},x_{1}^{-},\cdots,x_{n}^{+},x_{n}^{-}\}|\leq|\beta\cap\{x_{1}^{+},x_{1}^{-},\cdots,x_{n}^{+},x_{n}^{-}\}|$. 
In the case that $x$ is an outermost crossing of $D$, it follows that the image of $\alpha$ in $D'$ is a simple loop. 
We set $m=|\alpha\cap\{x_{1}^{+},x_{1}^{-},\cdots,x_{n}^{+},x_{n}^{-}\}|$. 
Then $m\leq n$. 
In the following, we consider the case that 
$|\alpha\cap\{x_{1}^{-},\cdots,x_{n}^{-}\}|\leq|\alpha\cap\{x_{1}^{+},\cdots,x_{n}^{+}\}|$. 
The other case is shown similarly. 
Then $\displaystyle{|\alpha\cap\{x_{1}^{-},\cdots,x_{n}^{-}\}|\leq\frac{m}{2}\leq\frac{n}{2}=\frac{c(D)-1}{2}}$. 
Let ${\mathbb S}^{2}$ be an equatorial $2$-sphere of ${\mathbb S}^{3}$ on which the diagram $D'$ is drawn. 
Let $U$ be an immersed circle in ${\mathbb S}^{2}$ obtained from $D'$ by forgetting over/under crossing information. 
By deforming $K'$ up to ambient isotopy if necessary, we may suppose that
\[
{\rm cl}(K'\setminus{\mathbb S}^{2})=\bigcup_{i=0}^{n}P_{i}
\]
and
\[
K'\subset U\cup\bigcup_{i=0}^{n}P_{i}
\]
where $P_{i}$ is an overpass contained in a small neighbourhood of $x_{i}$ with $x_{i}^{+}\in P_{i}$ for each $i\in\{0,1,\cdots,n\}$. 
The overpasses are contained in the upper hemisphere of ${\mathbb S}^{3}$ bounded by ${\mathbb S}^{2}$, and the points $x_{0}^{-},x_{1}^{-},\cdots,x_{n}^{-}$ are contained in $U$. 
We give an orientation to $K'$ such that the induced orientation on $\alpha$ is from $x_{0}^{+}$ to $x_{0}^{-}$. 
Let $p_{1},\cdots,p_{m}$ be the renaming of the points in $\alpha\cap\{x_{1}^{+},x_{1}^{-},\cdots,x_{n}^{+},x_{n}^{-}\}$ that appear in this order along the orientation of $\alpha$. 
Let $q_{0},r_{0},q_{1},r_{1},\cdots,q_{m},r_{m}$ be points in $\alpha$ such that the points $x_{0}^{+},q_{0},r_{0},p_{1},q_{1},r_{1},p_{2},\cdots,p_{m},q_{m},r_{m},x_{0}^{-}$ are arranged in this order along the orientation of $\alpha$. 
Let ${\mathcal O}$ (resp. ${\mathcal U}$) be the subset of $\{1,\cdots,m\}$ such that $p_{i}$ is an over-crossing (resp. under-crossing) if and only if $i\in{\mathcal O}$ (resp. $i\in{\mathcal U}$). 
By the assumption above we have $|{\mathcal U}|\leq\displaystyle{\frac{m}{2}}$. 
Let ${\mathcal P}$ (resp. ${\mathcal N}$) be the subset of $\{1,\cdots,m\}$ such that $p_{i}$ is a preimage of a positive crossing (resp. negative crossing) if and only if $i\in{\mathcal P}$ (resp. $i\in{\mathcal N}$). 
For $i\in\{1,\cdots,m\}$, let
\[
{\rm pn}(p_{i})=
\left\{ \begin{array}{ll}
1 & (i\in{\mathcal P})\\
-1 & (i\in{\mathcal N}).\\
\end{array} \right.
\]
Let ${\mathcal S}$ (resp. ${\mathcal M}$) be a subset of $\{1,\cdots,m\}$ such that $p_{i}$ is a preimage of a self crossing of $\alpha$ (resp. mutual crossing of $\alpha$ and $\beta$) if and only if $i\in{\mathcal S}$ (resp. $i\in{\mathcal M}$). 
We note that $K$ is ambient isotopic to a band sum of a Hopf link and $K'$ where the Hopf link and bands are contained in a small neighbourhood of $x_{0}$ in ${\mathbb S}^{3}$ as illustrated in Figure \ref{Hopf-band}. 
Then we lift up the Hopf link and bands so that they form an overbridge. 
Two vertical parts of the bands are denoted by $b_{1}$ and $b_{2}$ as illustrated in Figure \ref{Hopf-band2}. 
The band $b_{1}$ is attached to $K'$ at the point $q_{0}$ and the band $b_{2}$ is attached to $K'$ at the point $r_{m}$. 
We denote this form of $K$ by $K(q_{0},r_{m})$. 
For each $i\in\{0,1,\cdots,m\}$, we fix $i$ and consider the following deformation. 
We slide $b_{1}$, keeping it vertical, along the orientation of $\alpha$ to $q_{i}$. When $b_{1}$ encounter an over-crossing, $b_{1}$ passes through it by an ambient isotopy. When $b_{1}$ encounter an under-crossing, we perform a \lq\lq clasp leaps over a hurdle move\rq\rq \ illustrated in Figure \ref{hurdle} so that $b_{1}$ passes through the under-crossing. 
This move is realized by an application of delta move \cite[Claim 1.1]{M-N}\cite[Lemma 2.2]{T-Y}. 
In fact these two moves are equivalent as local move \cite[Examples]{T-Y2}. 
The horizontal parts of bands and a Hopf link are deformed up to ambient isotopy following the move of $b_{1}$. 
We denote the knots corresponding to the deformation so far by $K(q_{0},r_{m}),K(q_{1},r_{m}),\cdots,K(q_{i},r_{m})$. 
After $b_{1}$ is moved to $q_{i}$, we slide $b_{2}$, keeping it vertical, along the opposite orientation of $\alpha$ to $r_{i}$, up to ambient isotopy and applications of delta move just as above. 
We denote the knots corresponding to this deformation starting from $K(q_{i},r_{m})$ by $K(q_{i},r_{m-1}),K(q_{i},r_{m-2}),\cdots,K(q_{i},r_{i})$. 
The total number of applications of delta move so far is $|\alpha\cap\{x_{1}^{-},\cdots,x_{n}^{-}\}|=|{\mathcal U}|$. 
Namely we have $K(q_{i},r_{i})$ from $K=K(q_{0},r_{m})$ by $|{\mathcal U}|$-times application of delta move. 
Since $|{\mathcal U}|\leq\displaystyle{\frac{m}{2}\leq\frac{n}{2}=\frac{c(D)-1}{2}}$, we have $K(q_{i},r_{i})$ from $K$ by at most $\displaystyle{\left\lfloor\frac{c(D)-1}{2}\right\rfloor}$-times applications of delta move. 
In $K(q_{i},r_{i})$, a Hopf link and bands form a twist knot or a trivial knot as a factor knot connected summed to $K'$ at $q_{i}$ and $r_{i}$. 
Namely $K(q_{i},r_{i})$ is a connected sum of a twist knot or a trivial knot and $K'$. 
The knot type of this factor knot is determined by the total number of full-twists on the bands. 
This number is equal to the linking number of a $2$-component link obtained by smoothing any one of the two crossings of a Hopf link. 
In general, we denote the linking number of a $2$-component link obtained from the knot $K(q_{s},r_{t})$ by smoothing any one of the two crossings of a Hopf link by $l(s,t)$. 
In the following we will show
\begin{equation}
\label{lk1}
l(0,0)=-l(m,m) 
\end{equation}
and 
\begin{equation}
\label{lk2}
l(i+1,i+1)-l(i,i)\in\{-2,0,2\}
\end{equation}
for each $i\in\{0,1,\cdots,m-1\}$. 
Then by the intermediate value theorem, there exists $i\in\{0,1,\cdots,m\}$ such that $l(i,i)\in\{0,1\}$. 
Then the factor knot is $0_{1}$ or ${3_{1}}^*$. 
Suppose that the factor knot is $0_{1}$. Then $K(q_{i},r_{i})=K'$ and we are done. 
Suppose that the factor knot is ${3_{1}}^*$. Then $K(q_{i},r_{i})$ is a connected sum of ${3_{1}}^*$ and $K'$. 
Since $u_{\Delta}({3_{1}}^*)=1$, we have $K'$ from $K(q_{i},r_{i})$ by an application of delta move. 
Therefore we have $K'$ from $K$ by at most $\displaystyle{\left\lfloor\frac{c(D)-1}{2}\right\rfloor+1=\left\lfloor\frac{c(D)+1}{2}\right\rfloor}$-times applications of delta move. 
In case $x$ is outermost, we will show
\begin{equation}
\label{lk3}
l(0,0)=l(1,1)=\cdots=l(m,m)=0. 
\end{equation}
Then the factor knot is always $0_{1}$. Then $K(q_{0},r_{0})=K'$ and we are done. 
We consider the above stated sequence of $m+1$ knots
\[
K=K(q_{0},r_{m}),K(q_{1},r_{m}),\cdots,K(q_{i},r_{m}),K(q_{i},r_{m-1}),\cdots,K(q_{i},r_{i}).
\] 
Corresponding to this sequence of knots, we have a sequence of linking numbers 
\[
l(0,m),l(1,m),\cdots,l(i,m),l(i,m-1),\cdots,l(i,i).
\]
Let $d_{i,1},d_{i,2},\cdots,d_{i,m}$ be its difference sequence. 
The knots in the sequence above differ one by one, by the position of $b_{1}$ or $b_{2}$. 
For $j\leq i$, if $p_{j}$ is an over-crossing then $K(q_{j-1},r_{m})$ and $K(q_{j},r_{m})$ are ambient isotopic and therefore we have $d_{i,j}=l(j,m)-l(j-1,m)=0$. 
If $p_{j}$ is an under-crossing, then $K(q_{j},r_{m})$ is obtained from $K(q_{j-1},r_{m})$ by an application of delta move. 
Then we see $d_{i,j}\in\{-1,1\}$. This value $d_{i,j}$ depends on the value ${\rm pn}(p_{j})$ and the position of $b_{2}$ and ${\rm partner}(p_{j})$. 
In other words, it depends how the $6$ end points in Figure \ref{hurdle} are connected on the outside. 
Similarly, for $j>i$, if $p_{j}$ is an over-crossing then $K(q_{i},r_{j})$ and $K(q_{i},r_{j-1})$ are ambient isotopic and therefore $d_{i,i+1+m-j}=l(i,j-1)-l(i,j)=0$. 
If $p_{j}$ is an under-crossing, then $K(q_{i},r_{j-1})$ is obtained from $K(q_{i},r_{j})$ by an application of delta move and $d_{i,i+1+m-j}\in\{-1,1\}$. 
This value $d_{i,i+1+m-j}$ depends on the value ${\rm pn}(p_{j})$ and the position of $b_{1}$ and ${\rm partner}(p_{j})$. 
Let $\rho_{i}:\{1,\cdots,m\}\to\{1,\cdots,m\}$ be a bijection defined by
\[
\rho_{i}(j)=
\left\{ \begin{array}{ll}
j & (j\leq i)\\
i+1+m-j & (j>i).\\
\end{array} \right.
\]
Then we have
\[
l(i,i)=l(0,m)+\sum_{j=1}^{m}d_{i,j}=l(0,m)+\sum_{j\in{\mathcal U}}d_{i,\rho_{i}(j)}.
\]
By definition of linking number we have
\[
l(0,m)=\sum_{j\in{\mathcal U}\cap{\mathcal M}}{\rm pn}(p_{j}).
\]
For $j\in{\mathcal U}\cap{\mathcal M}$ we see, by observing how the $6$ end points in Figure \ref{hurdle} are connected on the outside, that $d_{i,\rho_{i}(j)}=-{\rm pn}(p_{j})$. 
Thus we have
\[
\sum_{j\in{\mathcal U}}d_{i,\rho_{i}(j)}=\sum_{j\in{\mathcal U}\cap{\mathcal M}}d_{i,\rho_{i}(j)}+\sum_{j\in{\mathcal U}\cap{\mathcal S}}d_{i,\rho_{i}(j)}=-l(0,m)+\sum_{j\in{\mathcal U}\cap{\mathcal S}}d_{i,\rho_{i}(j)}.
\]
Therefore we have
\[
l(i,i)=\sum_{j\in{\mathcal U}\cap{\mathcal S}}d_{i,\rho_{i}(j)}.
\]
In case $x$ is outermost, ${\mathcal S}$ is an empty set. Thus we have shown (\ref{lk3}). 
Let $j$ be an element of ${\mathcal U}\cap{\mathcal S}$. 
By observing how the $6$ end points in Figure \ref{hurdle} are connected on the outside, we see the following. 
Suppose that $j\neq i+1$ and ${\rm partner}(p_{j})\neq p_{i+1}$. 
Then we see $d_{i+1,\rho_{i+1}(j)}=d_{i,\rho_{i}(j)}$. 
Suppose that $j=i+1$ and ${\rm partner}(p_{j})=p_{k}$ with $k<j$. 
Then we see $d_{i+1,\rho_{i+1}(j)}=d_{i,\rho_{i}(j)}$. 
Suppose that $j=i+1$ and ${\rm partner}(p_{j})=p_{k}$ with $k>j$. 
Then we see $d_{i+1,\rho_{i+1}(j)}=-d_{i,\rho_{i}(j)}$. 
Suppose that ${\rm partner}(p_{j})=p_{i+1}$ and $j<i+1$. 
Then we see $d_{i+1,\rho_{i+1}(j)}=d_{i,\rho_{i}(j)}$. 
Suppose that ${\rm partner}(p_{j})=p_{i+1}$ and $j>i+1$. 
Then we see $d_{i+1,\rho_{i+1}(j)}=-d_{i,\rho_{i}(j)}$. 
These results imply that $d_{0,\rho_{0}(j)}=-d_{m,\rho_{m}(j)}$. 
Then we have
\[
l(0,0)=\sum_{j\in{\mathcal U}\cap{\mathcal S}}d_{0,\rho_{0}(j)}=\sum_{j\in{\mathcal U}\cap{\mathcal S}}(-d_{m,\rho_{m}(j)})=-\sum_{j\in{\mathcal U}\cap{\mathcal S}}d_{m,\rho_{m}(j)}=-l(m,m).
\]
Thus we have shown (\ref{lk1}). 
Since ${\rm partner}(p_{i+1})\neq p_{i+1}$, $j=i+1$ and ${\rm partner}(p_{j})=p_{i+1}$ do not occur simultaneously. 
Therefore we have shown (\ref{lk2}). 
This completes the proof. 
$\Box$


\vskip 5mm

Let $D$ be a knot diagram. We denote the minimal number of crossings of $D$ changing them turns $D$ to a diagram of a trivial knot by $u(D)$. 
It is well-known that $\displaystyle{u(D)\leq\frac{c(D)-1}{2}}$ if $c(D)>0$. See for example \cite{Taniyama}\cite[Proposition 2.1]{Ozawa2}. 
Moreover we have the following lemma. 

\begin{Lemma}\label{lemma-outermost-crossing} Let $D$ be a knot diagram. 
Suppose that there is a crossing of $D$ that is not outermost. 
Then
\[
u(D)\leq\frac{c(D)-3}{2}.
\]
\end{Lemma}

\noindent{\bf Proof.} Let $x_{0}$ be a crossing of $D$ that is not outermost. Let $C_{1}$ and $C_{2}$ be closed curves obtained from $D$ by smoothing $x_{0}$. Then both $C_{1}$ and $C_{2}$ are not simple closed curves. Let $x_{1}$ and $x_{2}$ be crossings of $C_{1}$ and $C_{2}$ respectively. 
Let $U$ be the underlying projection of $D$. Then we see that in the chord diagram of $U$, the chords corresponding to $x_{0}$, $x_{1}$ and $x_{2}$ are mutually parallel. Then we see that the trivializing number $t(U)$ of $U$ defined in \cite{Hanaki1} is less than or equal to $c(D)-3$ \cite[Lemma 3.7]{Henrich-etc}. It is implicitly shown in the proof of \cite[Proposition 1]{Hanaki2} that $\displaystyle{u(D)\leq\frac{t(U)}{2}}$. 
Therefore we have $\displaystyle{u(D)\leq\frac{c(D)-3}{2}}$. 
This completes the proof. 
$\Box$

\vskip 5mm

\noindent{\bf Proof of Theorem \ref{theorem-crossing-number-delta-unknotting-number}.} Let $K$ be a knot with $c(K)=4$. Then we have $K=4_{1}$ and $u_{\Delta}(K)=1$. Then
\[
u_{\Delta}(K)\leq\frac{1}{4}(c(K)^{2}-2c(K)-3)
\]
holds. 
Let $K$ be a knot with $c(K)\geq5$. 
First suppose that $K=T(2,p)$ for some odd number $p$ with $|p|\geq5$. 
It is known that $u_{\Delta}(T(2,p))=\displaystyle{\frac{p^{2}-1}{8}}$ \cite[Corollary 2.2]{N-N-U}. 
Since $c(T(2,p))=|p|$ it is sufficient to show
\[
\frac{|p|^{2}-1}{8}\leq\frac{|p|^{2}-2|p|-3}{4}
\]
for $|p|\geq5$. 
It is equivalent to 
\[
|p|^{2}-1\leq2|p|^{2}-4|p|-6. 
\]
Then it is equivalent to 
\[
|p|^{2}-4|p|-5\geq0.
\]
Since $|p|^{2}-4|p|-5=(|p|+1)(|p|-5)$ we have shown it. 
Next suppose that $K$ is not $T(2,p)$ for any odd number $p$. 
Let $D$ a diagram of $K$ with $c(D)=c(K)$. 
Then by \cite[Theorem 1.4 (1)]{Taniyama} we have $u(D)\leq\displaystyle{\frac{c(D)-2}{2}}$. 
Suppose that there is a crossing of $D$ that is not outermost. 
Then we have $\displaystyle{u(D)\leq\frac{c(D)-3}{2}}$ by Lemma \ref{lemma-outermost-crossing}. 
Then by Lemma \ref{lemma-crossing-change-by-delta-moves} we have
\[
u_{\Delta}(K)\leq(\frac{c(D)-3}{2})(\frac{c(D)+1}{2})=\frac{1}{4}(c(D)^{2}-2c(D)-3)=\frac{1}{4}(c(K)^{2}-2c(K)-3).
\]
Suppose that every crossing of $D$ is outermost. 
Then by Lemma \ref{lemma-crossing-change-by-delta-moves} we have
\[
u_{\Delta}(K)\leq(\frac{c(D)-2}{2})(\frac{c(D)-1}{2})=\frac{1}{4}(c(D)^{2}-3c(D)+2)=\frac{1}{4}(c(K)^{2}-3c(K)+2).
\]
Since $c(K)\geq5$ we have
\[
\frac{1}{4}(c(K)^{2}-3c(K)+2)\leq\frac{1}{4}(c(K)^{2}-2c(K)-3).
\]
This completes the proof. 
$\Box$

\vskip 5mm

\section{Braid index minus one and bridge number minus one}\label{section-bb} 

\noindent{\bf Proof of Theorem \ref{theorem-braid-index-bridge-number}.} Since ${\rm bridge}(K)\leq{\rm braid}(K)$ for every oriented knot $K$, $({\rm braid}-1,{\rm bridge}-1)(0_{1})=(0,0)$ and $({\rm bridge}-1)(K)=0$ if and only if $K$ is a trivial knot, we have
\[
({\rm braid}-1,{\rm bridge}-1)({\mathcal K})\subset\{(0,0)\}\cup\{(x,y)\in({\mathbb Z}_{>0})^{2}\mid y\leq x\}.
\]
As stated in Section \ref{section-cbraid}, we have, for a positive integer $a$, 
\[
({\rm braid}-1)(C(2,2a-1))=a.
\]
Since $C(2,2a-1)$ is a $2$-bridge knot we have
\[
({\rm bridge}-1)(C(2,2a-1))=1.
\]
Therefore we have, for a non-negative integer $b$, 
\[
({\rm braid}-1,{\rm bridge}-1)(C(2,2a-1)\#b\cdot{3_{1}}^*)=(a+b,1+b).
\]
Thus we have
\[
({\rm braid}-1,{\rm bridge}-1)({\mathcal K})\supset\{(0,0)\}\cup\{(x,y)\in({\mathbb Z}_{>0})^{2}\mid y\leq x\}.
\]
This completes the proof. 
$\Box$

\vskip 5mm

\section{Canonical genus and genus}\label{section-cgg} 

\noindent{\bf Proof of Theorem \ref{theorem-canonical-genus-genus}.} As stated in Section \ref{introduction}, it is shown in \cite[Theorem 1.1]{Stoimenow} that no knot has $g_{c}=2$ and $g=1$. Since $g(K)\leq g_{c}(K)$ for every oriented knot $K$, $(g_{c},g)(0_{1})=(0,0)$ and $g(K)=0$ if and only if $K$ is a trivial knot, we have 
\[
(g_{c},g)({\mathcal K})\subset\{(0,0)\}\cup(\{(x,y)\in({\mathbb Z}_{>0})^{2}\mid y\leq x\}\setminus\{(2,1)\}).
\]
Let $a$ be an integer with $a\geq3$. 
Let $K$ be a $2$-bridge knot with $c(K)=a$ and $D(K)$ a Whitehead double of $K$. 
Then $c_{g}(D(K))=c(K)=a$ \cite[Theorem 1]{Nakamura}. 
Let $\delta(J)$ be the maximal degree in $z$ of the HOMFLY-PT polynomial $P_{J}(v,z)$ of a knot $J$. 
It is shown in \cite{Morton} that $\displaystyle{g_{c}(J)\geq\frac{\delta(J)}{2}}$. 
It is also shown in \cite{Nakamura} that $\delta(D(K))=2c(K)=2a$. 
Since $\delta$ is additive under connected sum of knots and $\delta({3_{1}}^*)=2$, we have, for a non-negative integer $b$, 
\[
\delta(D(K)\#b\cdot{3_{1}}^*)=2c(K)+2b=2(a+b). 
\]
Therefore we have
\[
g_{c}(D(K)\#b\cdot{3_{1}}^*)=a+b. 
\]
Since genus is additive under connected sum of knots we have
\[
g(D(K)\#b\cdot{3_{1}}^*)=1+b. 
\]
Thus we have
\[
(g_{c},g)(D(K)\#b\cdot{3_{1}}^*)=(a+b,1+b). 
\]
Let ${\rm K11n39}$ be an $11$-crossing knot in Hoste-Thistlethwaite's table of 11 Crossing Knots. 
It is known that $g_{c}({\rm K11n39})=3$, $\delta({\rm K11n39})=6$ and $g({\rm K11n39})=2$. 
Therefore we have, for a non-negative integer $b$, 
\[
(g_{c},g)({\rm K11n39}\#b\cdot{3_{1}}^*)=(3+b,2+b). 
\]
We also have, for a non-negative integer $b$, 
\[
(g_{c},g)(b\cdot{3_{1}}^*)=(b,b). 
\]
As a summary of these results we have
\[
(g_{c},g)({\mathcal K})\supset\{(0,0)\}\cup(\{(x,y)\in({\mathbb Z}_{>0})^{2}\mid y\leq x\}\setminus\{(2,1)\}).
\]
This completes the proof. 
$\Box$

\vskip 5mm

\section{Independent pairs}\label{section-independent} 

\begin{Lemma}\label{lemma-trefoil-connected-sum} Let $a$ be a positive integer and $K$ an oriented knot in ${\mathbb S}^{3}$. 
Suppose $\sigma(K)=-2u(K)$. 
Then
\[
u(K\#a\cdot{3_{1}}^*)=u(K)+a.
\]
\end{Lemma}

\noindent{\bf Proof.} Since $u({3_{1}}^{*})=1$ we have $u(K\#a\cdot{3_{1}}^{*})\leq u(K)+a\cdot u({3_{1}}^{*})=u(K)+a$. 
Since $\sigma({3_{1}}^{*})=-2$ we have $\sigma(K\#a\cdot{3_{1}}^*)=\sigma(K)+a\cdot\sigma({3_{1}}^{*})=-2(u(K)+a)$. Therefore $u(K\#a\cdot{3_{1}}^{*})\geq u(K)+a$. Thus we have $u(K\#a\cdot{3_{1}}^{*})=u(K)+a$. 
This completes the proof. 
$\Box$

\vskip 5mm

In the following we consider local maximum and local minimum of a knot in ${\mathbb S}^{3}$. 
We also consider local maximum, local minimum and saddle of a torus in ${\mathbb S}^{3}$. 
They are defined with respect to the height function $h:{\mathbb S}^{3}\to{\mathbb R}$ sending $(x_{1},x_{2},x_{3},x_{4})\in{\mathbb S}^{3}$ to $x_{4}\in{\mathbb R}$. 
\begin{Lemma}\label{lemma-odd-bridge-knot} Let $P$ be a knot in ${\mathbb S}^{3}$ with ${\rm bridge}(P)\geq3$. 
Let $W$ be an unknotted solid torus in ${\mathbb S}^{3}$ with $P\subset{\rm int}W$. 
Suppose that $P$ is homotopically trivial in $W$ and the wrapping number of $P$ in $W$ is $2$. 
Let $X$ be a universal covering space of $W$ and $\psi:X\to W$ a universal covering projection. 
Suppose that each component of $\psi^{-1}(P)$ bounds a disk in ${\rm int}X$. 
Let $J$ be a non-trivial knot in ${\mathbb S}^{3}$. 
Let $K$ be a satellite knot in ${\mathbb S}^{3}$ with companion knot $J$ and pattern $(W,P)$. 
Then we have ${\rm bridge}(K)\geq2\cdot{\rm bridge}(J)+1$. 
\end{Lemma}

\noindent{\bf Proof.} 
By Schubert's theorem \cite{Schubert} re-proved in \cite{Schultens} we have 
$
{\rm bridge}(K)\geq2\cdot{\rm bridge}(J). 
$
We will show ${\rm bridge}(K)\geq2\cdot{\rm bridge}(J)+1$ by contradiction. 
Suppose that ${\rm bridge}(K)=2\cdot{\rm bridge}(J)$. 
We set $n={\rm bridge}(J)$ for simplicity. 
Let $V$ be a regular neighbourhood of $J$ in ${\mathbb S}^{3}$. 
Let $\varphi:W\to V$ be a faithful homeomorphism such that $K=\varphi(P)$. 
We note that the whole argument in \cite{Schultens} is applicable and we have the situation described below. 
We may suppose that $J$ is in minimal bridge-position. That is, $J$ has exactly $n$ local maximums and exactly $n$ local minimums. 
We may further suppose that $V$ is sufficiently thin so that the torus $T=\partial V$ has exactly $n$ local maximums, $n$ local minimums and $2n$ saddles. 
Each local maximum (resp. minimum) of $J$ has a small neighbourhood containing a local maximum (resp. minimum) of $T$ and a saddle of $T$. 
A meridian disk $D$ of $V$ is said to be {\it level} if $h(D)$ is a singleton. 
We choose mutually disjoint $4n$ level meridian disks $D_{1},\cdots,D_{4n}$ of $V$ and mutually interior-disjoint $4n$ $3$-balls $B_{1},\cdots,B_{4n}$ so that they satisfy the following conditions. 
Here we consider suffixes modulo $4n$. Namely we consider $4n+1=1$. 

\noindent
(1) $V=B_{1}\cup\cdots\cup B_{4n}$, 

\noindent
(2) $B_{i}\cap B_{i+1}=D_{i+1}$ ($i=1,\cdots,4n$), 

\noindent
(3) $B_{4j+1}$ contains exactly one local maximum of $J$ ($j=0,\cdots,n-1$), 

\noindent
(4) $\partial B_{4j+1}\cap T$ is an annulus containing exactly one local maximum of $T$ and exactly one saddle of $T$ ($j=0,\cdots,n-1$), 

\noindent
(5) $B_{4j+3}$ contains exactly one local minimum of $J$ ($j=0,\cdots,n-1$), 

\noindent
(6) $\partial B_{4j+3}\cap T$ is an annulus containing exactly one local minimum of $T$ and exactly one saddle of $T$ ($j=0,\cdots,n-1$). 

\noindent
Actually we take them so that $h(D_{4j+1})=h(D_{4j+2})$ and it is slightly lower than a local maximum of $J$ contained in $B_{4j+1}$, and 
$h(D_{4j+3})=h(D_{4j+4})$ and it is slightly higher than a local minimum of $J$ contained in $B_{4j+3}$ ($j=0,\cdots,n-1$). 
Then we may also suppose that $K\subset{\rm int}V$ satisfies the following conditions. 

\noindent
(7) $K$ has exactly $2n$ local maximums, 

\noindent
(8) $K$ intersects $D_{i}$ transversally at $2$ points for each $i\in\{1,\cdots,4n\}$, 

\noindent
(9) $B_{4j+1}$ contains exactly two local maximum of $K$ ($j=0,\cdots,n-1$), 

\noindent
(10) $B_{4j+3}$ contains exactly two local minimum of $K$ ($j=0,\cdots,n-1$). 

\begin{figure}[htbp]
      \begin{center}
\scalebox{0.5}{\includegraphics*{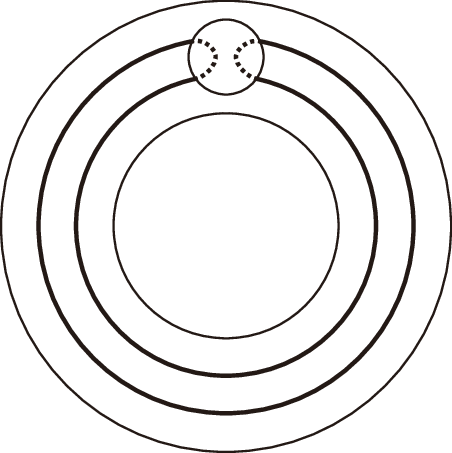}}
      \end{center}
   \caption{A pair $(V,K)$}
  \label{pattern2}
\end{figure}

\noindent
Then we see that both $B_{4j+2}$ and $B_{4j+4}$ contain no local maximums or local minimums of $K$ ($j=0,\cdots,n-1$). 
Then we see that $(B_{i},B_{i}\cap K)$ is a trivial tangle because it contains zero or two local maximums or local minimums of $K$ ($i=1,\cdots,4n$). 
We say that the tangle $(B_{i},B_{i}\cap K)$ is {\it forth-type} if each component of $B_{i}\cap K$ intersects $D_{i}$ at exactly one point. 
Otherwise it is said to be {\it back-type}. Since $K$ is connected and homotopically trivial in $V$, exactly one of $(B_{1},B_{1}\cap K),\cdots,(B_{4n},B_{4n}\cap K)$ is back-type. We see that both $(B_{4j+2},B_{4j+2}\cap K)$ and $(B_{4j+4},B_{4j+4}\cap K)$ are forth-type ($j=0,\cdots,n-1$). Then we may suppose without loss of generality that $(B_{1},B_{1}\cap K)$ is back-type. 
We now consider the universal covering projection $\varphi\circ\psi:X\to V$. A component of $(\varphi\circ\psi)^{-1}(K)$ bounds a disk in $X$. 
Since ${\rm int}X$ is homeomorphic to ${\mathbb R}^{3}$, it can be said that it is a trivial knot. 
By the situation described above, we may think that this trivial knot is a denominator of the tangle-sum of $2n-1$ tangles $(B_{3},B_{3}\cap K),(B_{5},B_{5}\cap K),\cdots,(B_{4n-1},B_{4n-1}\cap K)$. This knot is a connected sum of the denominators of these $2n-1$ tangles. 
This implies that each of the denominators of these trivial tangles must be a trivial knot. Then we see that these tangles are all integral tangles. 
Then, up to ambient isotopy of $V$, these tangles can be absorbed to the back-type tangle $(B_{1},B_{1}\cap K)$ and we see that $(V,K)$ is faithfully pairwise homeomorphic to the pair illustrated in Figure \ref{pattern2}. 
We note that the back-type tangle illustrated in Figure \ref{pattern2} is the result of absorption of integral tangles by the trivial tangle $(B_{1},B_{1}\cap K)$. Therefore it is still a trivial tangle. 
Since $(V,K)$ and $(W,P)$ are pairwise homeomorphic, we see that the knot $P$ is a numerator of this trivial tangle. Then we see that $P$ is a $2$-bridge knot. 
This contradicts to the assumption ${\rm bridge}(P)\geq3$. 
Thus we have ${\rm bridge}(K)\geq2n+1$. 
This completes the proof. 
$\Box$

\begin{Proposition}\label{proposition-unknotting-number-one-high-bridge-knot} Let $k$ be a positive integer. 
Then there exists an oriented knot $K$ with $(u,{\rm bridge}-1)(K)=(1,k)$ and $\sigma(K)=-2$. 
\end{Proposition}

\noindent{\bf Proof.} First suppose that $k$ is an odd number. 
For $k=1$ we have $(u,{\rm bridge}-1)({3_{1}}^{*})=(1,1)$ and $\sigma({3_{1}}^{*})=-2$. 
Suppose $k\geq3$. Then $k=2l+1$ for a positive integer $l$. 
Let $J=l\cdot{3_{1}}^{*}$ and $K$ a positively twisted double of $J$. Then we have $u(K)=1$. Since $K$ bounds a genus one Seifert surface whose Seifert matrix is same as that of a positive twist knot, we have $\sigma(K)=-2$. Since $({\rm bridge}-1)(J)=({\rm bridge}-1)(l\cdot{3_{1}}^{*})=l$ we have ${\rm bridge}(J)=l+1$. 
We may assume that $J$ is in minimal bridge-position. Then $J$ has exactly ${\rm bridge}(J)=l+1$ local maximums. 
By arranging $K$ in a thin regular neighbourhood of $J$ in ${\mathbb S}^{3}$ and by arranging the clasp of $K$ to a small neighbourhood of a local maximum of $J$, we have a representative of $K$ with $2(l+1)$ local maximums. Therefore ${\rm bridge}(K)\leq2(l+1)$. 
By Schubert's theorem we have ${\rm bridge}(K)\geq2(l+1)$. 
Thus we have ${\rm bridge}(K)=2(l+1)$. 
Therefore we have $(u,{\rm bridge}-1)(K)=(1,2l+1)=(1,k)$. 

\begin{figure}[htbp]
      \begin{center}
\scalebox{0.5}{\includegraphics*{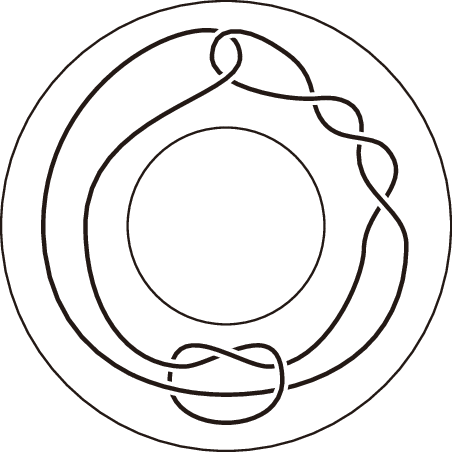}}
      \end{center}
   \caption{A pattern $P$ for unknotting number one odd-bridge knots}
  \label{pattern}
\end{figure}

\begin{figure}[htbp]
      \begin{center}
\scalebox{0.6}{\includegraphics*{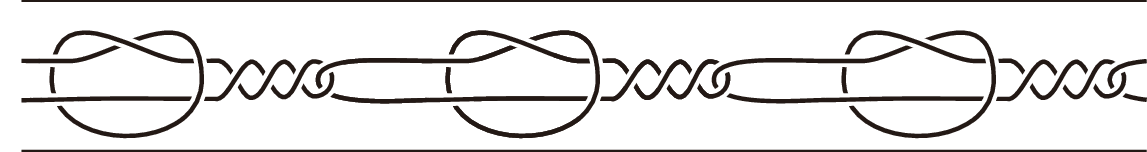}}
      \end{center}
   \caption{A lift of $P$}
  \label{lift}
\end{figure}

Next suppose that $k$ is an even number. 
For $k=2$ we have $(u,{\rm bridge}-1)({8_{21}}^{*})=(1,2)$ and $\sigma({8_{21}}^{*})=-2$. 
Suppose $k\geq4$. Then $k=2l+2$ for a positive integer $l$. 
Let $J=l\cdot3_{1}$ and $V$ a regular neighbourhood of $J$ in ${\mathbb S}^{3}$. 
Let $W$ be an unknotted solid torus in ${\mathbb S}^{3}$ and $P\subset{\rm int}W$ a pattern as illustrated in Figure \ref{pattern}. 
We see by a deformation up to ambient isotopy that $P={10_{133}}^{*}$ as a knot in ${\mathbb S}^{3}$. 
Then we have ${\rm bridge}(P)={\rm bridge}({10_{133}}^{*})=3$. 
Let $F\subset{\rm int}W$ be a Seifert surface of $P$ obtained by Seifert's algorithm applied to a diagram of $P$ illustrated in Figure \ref{pattern}. Let $\varphi:W\to V$ be a faithful homeomorphism. We set $K=\varphi(P)$. Namely $K$ is a satellite knot with companion knot $J$ and pattern $(W,P)$. 
Since the Seifert matrix of a Seifert surface $\varphi(F)$ of $K$ is same as that of $F$, we have $\sigma(K)=\sigma(P)=\sigma({10_{133}}^{*})=-2$. 
By removing the clasp of $K$ we have a trivial knot. Thus we have $u(K)=1$. 
We may assume that $J$ is in minimal bridge-position. Then $J$ has exactly ${\rm bridge}(J)=l+1$ local maximums. 
By arranging the clasp of $K$ to a small neighbourhood of a local maximum of $J$, and by arranging the $5$-crossing $2$-string tangle part of $P$ illustrated in Figure \ref{pattern} away from local maximums and local minimums of $J$, we have a representative of $K$ with $2(l+1)+1$ local maximums. Therefore ${\rm bridge}(K)\leq 2(l+1)+1$. 
Now we consider the universal covering space of $W$ and the lift of $P$ to it. They are as illustrated in Figure \ref{lift}.
We see that each component of the lift of $P$ is a trivial knot. 
Then by Lemma \ref{lemma-odd-bridge-knot} we have ${\rm bridge}(K)\geq2(l+1)+1$. Therefore we have $({\rm bridge}-1)(K)=2(l+1)=k$. 
This completes the proof. 
$\Box$


\vskip 5mm

\begin{Proposition}\label{proposition-genus-one-high-bridge-knot} Let $k$ be a positive integer. 
Then there exists an oriented knot $K$ with $(g,{\rm bridge}-1)(K)=(1,k)$. 
\end{Proposition}

\noindent{\bf Proof.} First suppose that $k$ is an odd number. 
For $k=1$ we have $(g,{\rm bridge}-1)({3_{1}}^{*})=(1,1)$. 
Suppose $k\geq3$. Then $k=2l+1$ for a positive integer $l$. 
Let $J=l\cdot{3_{1}}^{*}$ and $K$ a double of $J$. Then we have $g(K)=1$. 
We have shown in the proof of Proposition \ref{proposition-unknotting-number-one-high-bridge-knot} that ${\rm bridge}(K)=2l+2$. 
Therefore we have $(g,{\rm bridge}-1)(K)=(1,2l+1)=(1,k)$. 

\begin{figure}[htbp]
      \begin{center}
\scalebox{0.5}{\includegraphics*{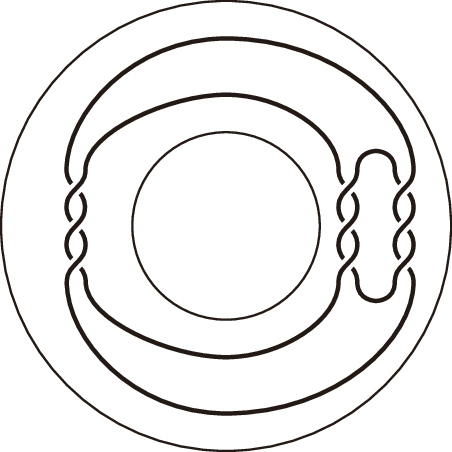}}
      \end{center}
   \caption{A pattern $Q$ for genus one odd-bridge knots}
  \label{pattern-genus}
\end{figure}

\begin{figure}[htbp]
      \begin{center}
\scalebox{0.6}{\includegraphics*{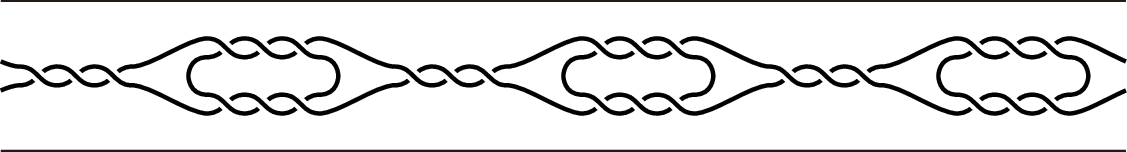}}
      \end{center}
   \caption{A lift of $Q$}
  \label{lift-genus}
\end{figure}

Next suppose that $k$ is an even number. 
For $k=2$ we have $(g,{\rm bridge}-1)(P(3,3,3))=(1,2)$. 
Suppose $k\geq4$. Then $k=2l+2$ for a positive integer $l$. 
Let $J=l\cdot3_{1}$ and $V$ a regular neighbourhood of $J$ in ${\mathbb S}^{3}$. 
Let $W$ be an unknotted solid torus in ${\mathbb S}^{3}$ and $Q\subset{\rm int}W$ a pattern as illustrated in Figure \ref{pattern-genus}. 
We see that $Q=P(3,3,3)$ as a knot in ${\mathbb S}^{3}$. 
Then we have ${\rm bridge}(Q)={\rm bridge}(P(3,3,3))=3$. 
Let $F\subset{\rm int}W$ be a Seifert surface of $P$ obtained by Seifert's algorithm applied to a diagram of $Q$ illustrated in Figure \ref{pattern-genus}. We note that $g(F)=1$. Let $\varphi:W\to V$ be a faithful homeomorphism. We set $K=\varphi(Q)$. Namely $K$ is a satellite knot with companion knot $J$ and pattern $(W,Q)$. Since $\varphi(F)$ is a Seifert surface of $K$ with $g(\varphi(F))=1$, we have $g(K)=1$. 
We may assume that $J$ is in minimal bridge-position. Then $J$ has exactly ${\rm bridge}(J)=l+1$ local maximums. 
We note that $Q$ has exactly $1$ local maximum in $W$ with respect to an ${\mathbb S}^{1}$ direction of $W$. 
Therefore we see ${\rm bridge}(K)\leq 2(l+1)+1$. 
Now we consider the universal covering space of $W$ and the lift of $Q$ to it. They are as illustrated in Figure \ref{lift-genus}.
We see that each component of the lift of $Q$ is a trivial knot. 
Then by Lemma \ref{lemma-odd-bridge-knot} we have ${\rm bridge}(K)\geq2(l+1)+1$. Therefore we have $(g,{\rm bridge}-1)(K)=(1,2(l+1))=(1,k)$. 
This completes the proof. 
$\Box$


\vskip 5mm

\noindent{\bf Proof of Theorem \ref{theorem-independent-pairs}.} First we show 
\[
(u,{\rm braid}-1)({\mathcal K})=\{(0,0)\}\cup({\mathbb Z}_{>0})^{2}.
\]
Since the conditions $u(K)=0$, $({\rm braid}-1)(K)=0$ and $K=0_{1}$ are mutually equivalent, we have
\[
(u,{\rm braid}-1)({\mathcal K})\subset\{(0,0)\}\cup({\mathbb Z}_{>0})^{2}.
\]
We note that $(u,{\rm braid}-1)(0_{1})=(0,0)$ and $(u,{\rm braid}-1)({3_{1}}^*)=(1,1)$. 
Let $k$ be a positive integer and $a$ a non-negative integer. 
Since $(u,{\rm braid}-1)(T(2,2k+1))=(k,1)$ and $\sigma(T(2,2k+1))=-2k=-2u(T(2,2k+1))$ we have by Lemma \ref{lemma-trefoil-connected-sum} 
\[
(u,{\rm braid}-1)(T(2,2k+1)\#a\cdot{3_{1}}^*)=(a+k,a+1). 
\]
Since $(u,{\rm braid}-1)(C(2,2k-1))=(1,k)$ and $\sigma(C(2,2k-1))=-2=-2u(C(2,2k-1))$ we have by Lemma \ref{lemma-trefoil-connected-sum} 
\[
(u,{\rm braid}-1)(C(2,2k-1)\#a\cdot{3_{1}}^*)=(a+1,a+k). 
\]
Thus we have
\[
(u,{\rm braid}-1)({\mathcal K})\supset\{(0,0)\}\cup({\mathbb Z}_{>0})^{2}
\]
and we have shown
\[
(u,{\rm braid}-1)({\mathcal K})=\{(0,0)\}\cup({\mathbb Z}_{>0})^{2}.
\]
%

Second we show
\[
(u,{\rm bridge}-1)({\mathcal K})=\{(0,0)\}\cup({\mathbb Z}_{>0})^{2}.
\]
Since the conditions $u(K)=0$, $({\rm bridge}-1)(K)=0$ and $K=0_{1}$ are mutually equivalent, we have
\[
(u,{\rm bridge}-1)({\mathcal K})\subset\{(0,0)\}\cup({\mathbb Z}_{>0})^{2}.
\]
We note that $(u,{\rm bridge}-1)(0_{1})=(0,0)$ and $(u,{\rm bridge}-1)({3_{1}}^*)=(1,1)$. 
Let $k$ be a positive integer and $a$ a non-negative integer. 
Since $(u,{\rm bridge}-1)(T(2,2k+1))=(k,1)$, 
we have 
\[
(u,{\rm bridge}-1)(T(2,2k+1)\#a\cdot{3_{1}}^*)=(a+k,a+1). 
\]
By Proposition \ref{proposition-unknotting-number-one-high-bridge-knot}, there exists an oriented knot $K$ with $(u,{\rm bridge}-1)(K)=(1,k)$ and $\sigma(K)=-2$. 
Then by Lemma \ref{lemma-trefoil-connected-sum} we have 
\[
(u,{\rm bridge}-1)(K\#a\cdot{3_{1}}^*)=(a+1,a+k). 
\]
Thus we have
\[
(u,{\rm bridge}-1)({\mathcal K})\supset\{(0,0)\}\cup({\mathbb Z}_{>0})^{2}
\]
and we have shown
\[
(u,{\rm bridge}-1)({\mathcal K})=\{(0,0)\}\cup({\mathbb Z}_{>0})^{2}.
\]
%

Third we show
\[
(g,{\rm braid}-1)({\mathcal K})=(g_{c},{\rm braid}-1)({\mathcal K})=\{(0,0)\}\cup({\mathbb Z}_{>0})^{2}.
\]
Since the conditions $g(K)=0$, $g_{c}(K)=0$, $({\rm braid}-1)(K)=0$ and $K=0_{1}$ are mutually equivalent, we have
\[
(g,{\rm braid}-1)({\mathcal K})\subset\{(0,0)\}\cup({\mathbb Z}_{>0})^{2}
\]
and
\[
(g_{c},{\rm braid}-1)({\mathcal K})\subset\{(0,0)\}\cup({\mathbb Z}_{>0})^{2}.
\]
Let 
\[
{\mathcal G}=\{0_{1}\}\cup\{C(2a,2b-1)\mid a,b\in{\mathbb Z}_{>0}\}.
\]
As we have mentioned in Section \ref{section-cuggc} and Section \ref{section-cbraid}, we have, for positive integers $a$ and $b$, 
\[
(g,{\rm braid}-1)(C(2a,2b-1))=(g_{c},{\rm braid}-1)(C(2a,2b-1))=(a,b).
\]
Therefore we have
\[
(g,{\rm braid}-1)({\mathcal K})\supset(g,{\rm braid}-1)({\mathcal G})\supset\{(0,0)\}\cup({\mathbb Z}_{>0})^{2}
\]
and
\[
(g_{c},{\rm braid}-1)({\mathcal K})\supset(g_{c},{\rm braid}-1)({\mathcal G})\supset\{(0,0)\}\cup({\mathbb Z}_{>0})^{2}
\]
and we have shown
\[
(g,{\rm braid}-1)({\mathcal K})=(g_{c},{\rm braid}-1)({\mathcal K})=\{(0,0)\}\cup({\mathbb Z}_{>0})^{2}.
\]
%

Last we show
\[
(g,{\rm bridge}-1)({\mathcal K})=\{(0,0)\}\cup({\mathbb Z}_{>0})^{2}.
\]
Since the conditions $g(K)=0$, $({\rm bridge}-1)(K)=0$ and $K=0_{1}$ are mutually equivalent, we have
\[
(g,{\rm bridge}-1)({\mathcal K})\subset\{(0,0)\}\cup({\mathbb Z}_{>0})^{2}.
\]
We note that $(g,{\rm bridge}-1)(0_{1})=(0,0)$ and $(g,{\rm bridge}-1)({3_{1}}^*)=(1,1)$. 
Let $k$ be a positive integer and $a$ a non-negative integer. 
Since $(g,{\rm bridge}-1)(T(2,2k+1))=(k,1)$, 
we have 
\[
(g,{\rm bridge}-1)(T(2,2k+1)\#a\cdot{3_{1}}^*)=(a+k,a+1). 
\]
By Proposition \ref{proposition-genus-one-high-bridge-knot}, there exists an oriented knot $K$ with $(g,{\rm bridge}-1)(K)=(1,k)$. 
Then we have 
\[
(g,{\rm bridge}-1)(K\#a\cdot{3_{1}}^*)=(a+1,a+k). 
\]
Thus we have
\[
(g,{\rm bridge}-1)({\mathcal K})\supset\{(0,0)\}\cup({\mathbb Z}_{>0})^{2}
\]
and we have shown
\[
(g,{\rm bridge}-1)({\mathcal K})=\{(0,0)\}\cup({\mathbb Z}_{>0})^{2}.
\]
This completes the proof. 
$\Box$

\vskip 5mm

\begin{figure}[htbp]
      \begin{center}
\scalebox{0.5}{\includegraphics*{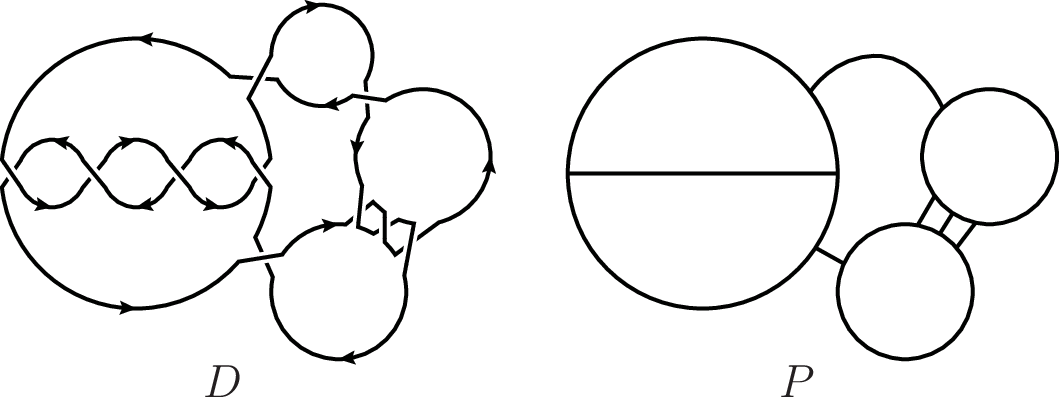}}
      \end{center}
   \caption{A $1$-complex $P$ for a knot diagram $D$}
  \label{complex}
\end{figure}

\noindent{\bf Proof of Theorem \ref{theorem-canonical-genus-bridge-number}.}
Let $k$ be a positive integer. 
Let $K$ be an oriented knot in ${\mathbb S}^{3}$ with $g_{c}(K)=k$. 
Then $K$ has a diagram $D$ on ${\mathbb S}^{2}$ with $g(F(D)))=k$. 
We may choose such $D$ so that it has no nugatory crossings. 
Then each nested Seifert circle of $D$ has at least $4$ crossings of $D$. 
Let ${\mathcal C}(D)$ be the set of all crossings of $D$ and ${\mathcal S}(D)$ the set of all Seifert circles of $D$. 
A Seifert graph $G(D)$ of $D$ is the graph with the vertex set ${\mathcal S}(D)$ and the edge set ${\mathcal C}(D)$. 
Each edge of $G(D)$, that is a crossing of $D$, joins two vertices of $G(D)$, that are Seifert circles of $D$ facing each other at that crossing. 
Let $V_{i}(D)$ be the number of degree $i$ vertices of $G(D)$. We note that $G(D)$ is a connected graph and $V_{0}(D)=0$. 
Since $D$ has no nugatory crossings, $V_{1}(D)=0$. 
Then by shake-hands lemma of graph theory we have 
\[
\sum_{i\geq2}i\cdot V_{i}(D)=2c(D).
\]
Since 
\[
\sum_{i\geq2}V_{i}(D)=s(D)
\]
we have 
\[
\chi(F(D))=s(D)-c(D)=\frac{1}{2}(2\sum_{i\geq2}V_{i}(D)-\sum_{i\geq2}i\cdot V_{i}(D)).
\]
Since
\[
2\sum_{i\geq2}V_{i}(D)-\sum_{i\geq2}i\cdot V_{i}(D)=\sum_{i\geq2}2V_{i}(D)-\sum_{i\geq2}i\cdot V_{i}(D)=\sum_{i\geq3}(2-i)V_{i}(D)
\]
and
\[
\chi(F(D))=1-2g(F(D))=1-2k
\]
we have 
\[
1-2k=\frac{1}{2}(\sum_{i\geq3}(2-i)V_{i}(D)).
\]
Therefore
\[
\sum_{i\geq3}V_{i}(D)\leq\sum_{i\geq3}(i-2)V_{i}(D)=4k-2.
\]
Thus we see that $D$ has at most $4k-2$ Seifert circles with $3$ or more crossings. 
As stated above, each Seifert circle with exactly $2$ crossings are not nested, and therefore innermost on ${\mathbb S}^{2}$. 
Then we can suppose that $D$ is contained in a thin regular neighbourhood of a $1$-complex $P$ in ${\mathbb S}^{2}$ that is obtained from a disjoint union of all Seifert circles of $D$ by adding simple arcs each of which corresponds to a crossing of $D$ connecting two Seifert circles facing each other at that crossing, and then contracting each Seifert circle with exactly $2$ crossings to a point. See for example Figure \ref{complex}. 
We see that $P$ is homeomorphic to a $3$-regular graph embedded in ${\mathbb S}^{2}$. 
Let $V(P)$ be the number of vertices of $P$ and $E(P)$ the number of edges of $P$. 
By shake-hands lemma we have $3V(P)=2E(P)$. 
Since $P$ is a deformation retract of $F(D)$ we have $\chi(P)=\chi(F(D))=1-2k$. 
On the other hand $\displaystyle{\chi(P)=V(P)-E(P)=-\frac{1}{2}V(P)}$. Then we have $V(P)=4k-2$. 
We note that there exists only finitely many $3$-regular graphs with $4k-2$ vertices embedded in ${\mathbb S}^{2}$ up to ambient isotopy of ${\mathbb S}^{2}$. 
Each of such $3$-regular graph provides a common upper bound for the bridge number of a knot represented by a diagram contained in a thin neighbourhood of the graph as described above. Let $m$ be the maximum of the upper bounds for the bridge number of all such $3$-regular graphs. 
Then every oriented knot $K$ with $g_{c}(K)=k$ satisfies ${\rm bridge}(K)\leq m$. 
This completes the proof. 
$\Box$

\section{Future directions}\label{section-future} 

Let $\alpha:{\mathcal K}\to X$ be an oriented knot invariant and ${\mathcal J}$ a subset of ${\mathcal K}$. 
Then $\alpha({\mathcal J})$ is a subset of $\alpha({\mathcal K})$. 
In some cases $\alpha({\mathcal J})=\alpha({\mathcal K})$ and in some cases $\alpha({\mathcal J})$ is a proper subset of $\alpha({\mathcal K})$. 
For example we have remarked in Section \ref{introduction} that 
\[
\Delta({\mathcal K})=\Delta(\{K\in{\mathcal K}\mid u(K)=1\}).
\]
We now extend this study to pairs of knot invariants. 
Let $\beta:{\mathcal K}\to Y$ be another oriented knot invariant. Then $(\alpha,\beta)({\mathcal J})$ is a subset of $(\alpha,\beta)({\mathcal K})$. 
The determination of the set
$(\alpha,\beta)({\mathcal K})\setminus(\alpha,\beta)({\mathcal J})$
may be an interesting problem. 
It is shown in Theorem \ref{theorem-crossing-number-unknotting-number-genus-canonical-genus} that $(c,u)({\mathcal K})\setminus(c,u)({\mathcal P}_{0})=\{(6,2)\}=(c,u)(\{{3_{1}}^*\#{3_{1}}^*\})$. It is also shown in Theorem \ref{theorem-crossing-number-unknotting-number-genus-canonical-genus} that both $(c,g)({\mathcal K})\setminus(c,g)({\mathcal R}_{0})$ and $(c,g_{c})({\mathcal K})\setminus(c,g_{c})({\mathcal R}_{0})$ are empty sets. 
For the sake of simplicity we have used connected sum construction of knots for proofs of some theorems in this paper. 
It may be possible to check that $(\alpha,\beta)({\mathcal K})\setminus(\alpha,\beta)({\mathcal P}_{0})$ is an empty set for some pairs $(\alpha,\beta)$. 

The concept pair of knot invariants is naturally generalized to triple of knot invariants, quadruple of knot invariants, and in general, $n$-tuple of knot invariants. 
Let $X_{i}$ be a set and $\alpha_{i}:{\mathcal K}\to X_{i}$ an oriented knot invariant for $i\in\{1,\cdots,n\}$. 
Let 
\[
(\alpha_{1},\cdots,\alpha_{n}):{\mathcal K}\to X_{1}\times\cdots\times X_{n}
\]
be a map defined by 
\[
(\alpha_{1},\cdots,\alpha_{n})(K)=(\alpha_{1}(K),\cdots,\alpha_{n}(K))
\]
for $K\in{\mathcal K}$. 
The subset $(\alpha_{1},\cdots,\alpha_{n})({\mathcal K})$ of $X_{1}\times\cdots\times X_{n}$ is said to be the relation of $(\alpha_{1},\cdots,\alpha_{n})$. 
The following example shows that the triple 
\[
(c,u,{\rm braid}-1):{\mathcal K}\to({\mathbb Z}_{\geq0})^{3}
\]
is not restored by the pairs 
\[
(c,u):{\mathcal K}\to({\mathbb Z}_{\geq0})^{2},
\]
 \[
 (c,{\rm braid}-1):{\mathcal K}\to({\mathbb Z}_{\geq0})^{2}
 \]
 and 
 \[
 (u,{\rm braid}-1):{\mathcal K}\to({\mathbb Z}_{\geq0})^{2}.
 \]
We have $(c,u)(5_1)=(5,2)$, $(c,{\rm braid}-1)(5_2)=(5,2)$ and $(u,{\rm braid}-1)({3_{1}}^*\#{3_{1}}^*)=(u,{\rm braid}-1)(7_{3})=(2,2)$. However there exist no oriented knot $K$ with $(c,u,{\rm braid}-1)(K)=(5,2,2)$. 

As an example of the relation of a triple of knot invariants, we note here that an inequality
\[
c(K)\geq2g(K)+({\rm braid}-1)(K)
\]
holds for every knot $K$ \cite[Theorem 2.6]{Diao}. The proof actually shows that a stronger inequality
\[
c(K)\geq2g_{c}(K)+({\rm braid}-1)(K)
\]
holds for every knot $K$. 
As a relevant fact, we note that the knots in Figure \ref{cg-knots} and the knots in Figure \ref{cbraid-1-knots} are exactly the same. They are just upside down each other. 
We also note that Quantitative Birman-Menasco finiteness theorem in \cite{Ito} after Birman-Menasco finiteness theorem in \cite{B-M2} gives an opposite estimation of crossing number by genus and braid index. 

In \cite{L-L} a result involving unknotting number, genus and braid index is shown. 

In \cite{Willerton} the pair of order-two Vassiliev invariant and order-three Vassiliev invariant of a fixed crossing number is exhibited. 
It is called Willerton's fish. It is a cross-section of the triple of crossing number, order-two Vassiliev invariant and order-three Vassiliev invariant. 

We see in Example \ref{example-genus-Euler-characteristic} that the genus and the Euler characteristic of closed connected orientable surfaces are mutually dependent and $(\chi,g)({\mathcal T})$ seems to be $1$-dimensional. See Figure \ref{genusEuler}. On the other hand, most of $(\alpha,\beta)({\mathcal K})$ in this paper seems $2$-dimensional. The following formulation verifies this intuition. 
Let $\alpha_{i}:{\mathcal K}\to{\mathbb Z}$ be an oriented knot invariant for $i\in\{1,\cdots,n\}$. 
Let $m$ be a non-negative integer. 
We say that the dimension ${\rm dim}(\alpha_{1},\cdots,\alpha_{n})$ is less than $m$ if 
\[
\lim_{k\to\infty}|(\alpha_{1},\cdots,\alpha_{n})({\mathcal K})\cap[-k,k]^{n}|/k^{m}=0.
\]
Suppose that ${\rm dim}(\alpha_{1},\cdots,\alpha_{n})$ is not less than $m$ but less than $m+1$. 
Then we define ${\rm dim}(\alpha_{1},\cdots,\alpha_{n})=m$. 
As an example we will show that 
\[
{\rm dim}(c,{\rm braid}-1,{\rm bridge}-1)=3.
\]
We note that $(c,{\rm braid}-1,{\rm bridge}-1)({3_{1}}^*)=(3,1,1)$,\\
$(c,{\rm braid}-1,{\rm bridge}-1)(4_{1})=(4,2,1)$ and $(c,{\rm braid}-1,{\rm bridge}-1)(5_{1})=(5,1,1)$. 
Then we have 
\[
(c,{\rm braid}-1,{\rm bridge}-1)(p\cdot{3_{1}}^*\#q\cdot4_{1}\#r\cdot5_{1})=p(3,1,1)+q(4,2,1)+r(5,1,1).
\] 
Since vectors $(3,1,1)$, $(4,2,1)$ and $(5,1,1)$ are linearly independent, we have 
\[
{\rm dim}(c,{\rm braid}-1,{\rm bridge}-1)=3.
\]

\vskip 3mm

\section*{Acknowledgments} The author is grateful to the referee for letting him know geography and botany problems in mathematics. 

\vskip 3mm

{\normalsize

}

\end{document}